\documentclass[aop,MSNbibl,seceqn,dvips]{arximspdf}
\usepackage{url,breakurl}

\doi{10.1214/13-AOP893} 
\volume{43}
\issue{1}
\pubyear{2015}
\firstpage{44}
\lastpage{84}

\makeatletter
\def\cal{\mathcal}
\def\mid{\vert}
\newproclaim{example}{Example}[section]
\newtheorem{theorem}{Theorem}[section]
\newtheorem{lemma}{Lemma}[section]
\newtheorem{corollary}{Corollary}[section]

\newproclaim{remark}{Remark}[section]
\newproclaim{inition}{Definition}[section]

\def\al{\alpha}
\def\bb{\beta}
\def\ga{\gamma}
\def\De{\Delta}
\def\ep{\varepsilon}
\def\la{\lambda}
\def\ka{\kappa}
\def\Om{\Omega}
\def\si{\sigma}
\def\th{\theta}
\def\ze{\zeta}

\def\rar{\rightarrow}
\def\cd{ \cdot}

\def\FF{{\cal F}}
\def\LL{{\cal L}}
\def\MM{{\cal M}}
\def\OO{{\cal O}}
\def\RR{{\cal R}}
\def\VV{{\cal V}}

\makeatother

\begin{document}
\begin{frontmatter}

\title{Permanental fields, loop soups and continuous additive functionals}
\runtitle{Permanental fields, loop soups and CAFs}

\begin{aug}
\author[A]{\fnms{Yves} \snm{Le~Jan}\ead[label=e1]{yves.lejan@math.u-psud.fr}},
\author[B]{\fnms{Michael B.} \snm{Marcus}\ead[label=e2]{mbmarcus@optonline.net}\thanksref{T1}}
\and
\author[C]{\fnms{Jay} \snm{Rosen}\corref{}\ead[label=e3]{jrosen30@optimum.net}\thanksref{T1,T2}}
\thankstext{T1}{Supported in part by grants from the National Science
Foundation and
PSCCUNY.}
\thankstext{T2}{Supported in part by a grant from the Simons Foundation.}
\runauthor{Y. Le Jan, M.~B. Marcus and J. Rosen}
\affiliation{Universit\'e Paris-Sud, City University of New York and\break
City University of New York}
\address[A]{Y. Le Jan\\
Equipe Probabilit\'es et Statistiques\\
Universit\'e Paris-Sud\\
B\^atiment 425\\
91405 Orsay Cedex\\
France\\
\printead{e1}}
\address[B]{M.~B. Marcus\\
Department of Mathematics\\
City College\\
City University of New York\\
New York, New York 10031\\
USA\\
\printead{e2}}
\address[C]{J. Rosen\\
Department of Mathematical\\
College of Staten Island\\
City University of New York\\
Staten Island, New York 10314\\
USA\\
\printead{e3}}
\end{aug}

\received{\smonth{9} \syear{2012}}
\revised{\smonth{8} \syear{2013}}

%
\begin{abstract}

A permanental field, $\psi=\{\psi(\nu), \nu\in\VV\}$, is a particular
stochastic process indexed by a space of measures on a set $S$. It is
determined by a kernel $u(x,y)$, $x,y\in S$, that need not be symmetric
and is allowed to be infinite on the diagonal. We show that these
fields exist when $u(x,y)$ is a potential density of a transient Markov
process $X$ in $S$.

A permanental field $\psi$ can be realized as the limit of a
renormalized sum of continuous additive functionals determined by a
loop soup of $X$, which we carefully construct. A Dynkin-type
isomorphism theorem is obtained that relates $\psi$ to continuous
additive functionals of $X$ (continuous in~$t$), $L=\{L_{t}^{\nu},
(\nu
,t)\in\VV\times R_{+}\}$. Sufficient conditions are obtained for the
continuity of $L$ on $\VV\times R_{+}$. The metric on $\VV$ is given by
a \textit{proper norm}.
\end{abstract}

%


\end{frontmatter}

\section{Introduction}\label{sec1}

In \cite{MR96}, we use a version of the Dynkin isomorphism theorem to
analyze families of continuous additive functionals of symmetric Markov
processes in terms of associated second-order Gaussian chaoses that are
constructed from Gaussian fields with covariance kernels that are the
potential densities of the symmetric Markov processes.

In this paper, we define a permanental field, $\psi=\{\psi(\nu), \nu
\in
\VV\}$, a new stochastic process indexed by a space of measures $\VV$
on a set $S$, that is determined by a kernel $u(x,y)$, $x,y\in S$, that
need not be symmetric. Permanental fields are a generalization of
second-order Gaussian chaoses. We show that these fields exist whenever
$u(x,y)$ is the potential density of a transient Markov process $X$.

We show that $\psi$ can be realized as the limit of a renormalized sum
of continuous additive functionals determined by a loop soup of $X$. A
loop soup is a Poisson point process on the path space of $X$ with an
intensity measure $\mu$ called the ``loop measure.'' (This is done in
Section~\ref{sec-exist}.) We obtain a new Dynkin type isomorphism
theorem that relates $\psi$ to continuous additive functionals of $X$
and can be used to analyze them.

Let $(\Om, \mathcal{F},P)$ be a probability space, and let $
S$ be a locally compact metric space with countable base. Let $\mathcal
{B}( S)$ denote the Borel $\si$-algebra, and let $\mathcal{M}( S)$ be
the set of finite signed Radon measures on $\mathcal{B}( S)$.

\begin{inition}\label{def1}
A map $\psi$ from a subset $\mathcal{V}\subseteq\mathcal{M}( S)$ to
$\mathcal{F}$ measurable functions on $\Om$ is called an $\al
$-permanental field with kernel $u$ if for all $\nu\in\VV$, $E\psi
(\nu
)=0$ and for all integers $n\ge2$ and $\nu_{1}, \ldots, \nu_{n} \in
\mathcal{V}$
%
\begin{equation}
E\Biggl( \prod_{j=1}^{n}\psi(
\nu_{j}) \Biggr) = \sum_{\pi\in\mathcal{P}'}\al
^{c(\pi)} \int\prod_{j=1}^{n}u
(x_{j},x_{\pi(j)})\prod_{j=1}^{n}
\,d\nu _{j}(x_{j}),\label{a.3}
\end{equation}
where $\mathcal{P}'$ is the set of permutations $\pi$ of $[1,n]$ such
that $ \pi(j)\neq j$ for any $j$, and $c(\pi)$ is the number of cycles
in the permutation $\pi$.
\end{inition}

The concept of permanental fields is motivated by \cite{LeJan}
and \cite{LeJan1}, Chapter~9.

The statement in (\ref{a.3}) makes sense when the kernel $u$
is bounded. However, in this case, we can accomplish the goals of this
paper using permanental processes as we do in
\cite{MRperm}. In this paper, we are particularly interested in the
case in which~$u$ is infinite on the diagonal. That is why we define
the field using measures on $S$ rather than points in $S$, and require
that $ \pi(j)\neq j$ for any $j$ in (\ref{a.3}) [since we allow
$u(x_{j},x_{j})=\infty$].

When $u$ is symmetric, positive definite and $\al=1/2$, $\{\psi(\nu),
\nu\in\VV\}$
is given by the Wick square, a particular second-order Gaussian chaos
defined as
%
\begin{equation}
:G^{2}: (\nu)=\lim_{\delta\rar0} \int\bigl(
G^{2}_{x,\delta
}-E\bigl(G^{2}_{x,\delta}\bigr)
\bigr) \,d\nu(x), \label{a.3a}
\end{equation}
where
$\{G_{x,\delta}, x\in S\} $ is a mean zero Gaussian process with finite
covariance $u_{\delta}(x,y)$, and $\lim_{\delta\rar0} u_{\delta
}(x,y)=u(x,y)$.
(See \cite{MR96} for details.) The results in \cite{MR96} are simpler
to achieve than the results in this paper because we have at our
disposal a wealth of information about second-order Gaussian chaoses.

The definition of a permanental field in (\ref{a.3}) is a
generalization of the moment formula for permanental processes,
introduced in \cite{VJ}. Let
$\th=\{\th_{x}, x\in S\} $ be an $\al$-permanental process with
(finite) kernel $u$, then for any $x_{1},\ldots, x_{n}\in S$
%
\begin{equation}
E\Biggl( \prod_{j=1}^{n}
\th_{x_{j}}\Biggr) = \sum_{\pi\in\mathcal{P}}\al
^{c(\pi)} \prod_{j=1}^{n}u
(x_{j},x_{\pi(j)}), \label{1.1m}
\end{equation}
where $\mathcal{P}$ is the set of permutations $\pi$ of $[1,n]$, and
$c(\pi)$ is the number of cycles in the permutation $\pi$. In this
case, $ \int(\th_{x}-E(\th_{x})) \,d\nu(x)$ is a permanental field.

Eisenbaum and Kaspi \cite{EK} show
that an $\al$-permanental process with kernel $u$ exists whenever $u$
is the potential density of a transient Markov process $X$ in $S$.
(This can also be done using loop soups. See \cite{LeJan1}, Chapters~2, 4, 5, for a study in the discrete symmetric case.) In \cite{MRperm}, we
give sufficient conditions for the continuity of $\al$-permanental
processes and use this, together with an isomorphism theorem of
Eisenbaum and Kaspi, \cite{EK} to give sufficient conditions for the
joint continuity of the local times of $X$. In this paper, we extend
these results to permanental fields and continuous additive functionals.

In order that (\ref{a.3}) makes sense, we need bounds on multiple
integrals of the form
%
\begin{equation}
\int\prod_{j=1}^{n}u (x_{j},x_{j+1})
\prod_{i=1}^{n} \,d\nu_{j}
(x_{i} ),\qquad x_{n+1}=x_{1}.
\end{equation}
We say that a norm $\|\cdot\|$ on $\MM(S)$ is a \textit{proper norm} with
respect to a kernel $u$ if for all $n\geq2$ and $\nu_{1},\ldots, \nu
_{n}$ in $\MM(S)$
%
\begin{equation}
\Biggl| \int\prod_{j=1}^{n}u
(x_{j},x_{j+1})\prod_{i=1}^{n}
\,d\nu_{j} (x_{i} ) \Biggr| \leq C^{n}\prod
_{j=1}^{n} \|\nu_{j}\|\label{p.0}
\end{equation}
for some universal constant $C<\infty$.

In Section~\ref{sec-levy}, in which we consider the continuity of
certain additive functionals of L\'evy processes, an explicit example
of a proper norm is given in (\ref{bp.1q}). Another example of a proper
norm which plays an important role in this paper is given in~(\ref
{p.2}). Additional examples of proper norms are given in Example~\ref
{exam-6.2}.

The next step in our program is to show that permanental fields exist.
We do this in Section~\ref{sec-exist} when the kernel $u(x,y)$ is the
potential density of a transient Borel right process $X$ in $S$.
(Additional technical conditions are given in Section~\ref{ss-caf}.)

We denote by $\mathcal{R}^{+}(X)$, or $\mathcal{R}^{+}$ when $X$ is
understood, the set of positive bounded Revuz measures $\nu$ on $S$
that are associated with $X$. This is explained in detail in Section~\ref{ss-caf}.

Let $\|\cdot\|$ be a proper norm on $\MM(S)$ with respect to the kernel
$u $. Set
%
\begin{equation}
\MM^{+}_{\|\cdot\|}=\bigl\{\mbox{positive }\nu\in\MM(S) | \|\nu\|
<\infty\bigr\}\label{con.1i}
\end{equation}
and
%
\begin{equation}
\mathcal{R}^{+}_{\|\cdot\|}=\mathcal{R}^{+}\cap
\MM^{+}_{\|\cdot\|
}.\label{rev1i}
\end{equation}
Let $\MM_{\|\cdot\|}$ and $\mathcal{R}_{\|\cdot\|}$ denote the set of
measures of the form $\nu=\nu_{1}-\nu_{2}$ with $ \nu_{1},\nu
_{2}\in
\MM^{+}_{\|\cdot\|}$ or $ \mathcal{R}^{+}_{\|\cdot\|}$, respectively.
We often omit saying that both $\RR_{\|\cdot\|}$ and $\|\cdot\|$ depend
on the kernel $u$.

The following theorem is implied by the results in Section~\ref{sec-exist}.

\begin{theorem}\label{theo-9.1i} Let $X$ be a transient Borel right
process with
state space $S$ and potential density $u(x,y)$, $x,y\in S$, as
described in Section~\ref{ss-caf}, and let $\|\cdot\|$ be a proper norm
with respect to the kernel $u(x,y)$. Then for $\al>0$ we can find an
$\al$-permanental field $\{\psi(\nu),\nu\in\mathcal{R}_{\|\cdot
\|}\}
$ with kernel $u$.
\end{theorem}

We say that $\{\psi(\nu),\nu\in\mathcal{R}_{\|\cdot\|}\}$ is the
$\al$-permanental field associated with~$X$.

In Section~\ref{sec-contperm} we study the continuity of
permanental fields. Let $\{\psi(\nu),\nu\in\VV\}$ be a permanental
field with kernel $u$. Let $\|\cdot\|$ be a proper norm with respect
to $u$ and suppose that $\VV\subseteq\MM_{\|\cdot\|}$. We show in
Section~\ref{sec-contperm} that
%
\begin{equation}
\bigl\| \psi(\mu)-\psi(\nu)\bigr\|_{\Xi}\le C \| \mu-\nu\|,\label{1.9}
\end{equation}
where $ \| \cd\|_{\Xi}$ is the norm of the exponential Orlicz space
generated by $e^{|x|}-1$.
This inequality enables us to use the well-known majorizing measure
sufficient condition for the continuity of stochastic processes, to
obtain sufficient conditions for the continuity of permanental fields,
$\{\psi(\nu),\nu\in\VV\}$ on $(\VV, \|\cdot\|)$, where $\|\cdot
\|$
denotes the metric $\|\mu-\nu\|$ in (\ref{1.9}).

Let $B_{\|\cdot\|}(\nu,r)$ denote the closed ball in $(\VV, \|\cdot
\|)$
with radius $r$ and center $\nu$. For any probability measure $\si$ on
$(\VV, \|\cdot\|)$, let
%
\begin{equation}
J_{\VV, \|\cdot\|,\si}( a) =\sup_{\nu\in\VV}\int_0^a
\log\frac{1}{\si
(B_{\|\cdot\|}(\nu,r))} \,dr.\label{taui}
\end{equation}

\begin{theorem}\label{theo-contprop} Let $\{\psi(\nu),\nu\in\VV\}
$ be an
$\al
$-permanental field with kernel $u$ and let $\|\cdot\|$ be a proper
norm for $u$. Assume that there exists
a probability measure $\si$ on $\VV$ such that $J_{\VV, \|\cdot\|
,\si
}(D)<\infty$, where $D$ is the diameter of $\VV$ with respect to $ \|
\cdot
\|$ and
%
\begin{equation}
\lim_{\delta\to0} J_{\VV, \|\cdot\|,\si}(\delta)=0. \label{2.10vix}
\end{equation}
Then $\psi$
is uniformly continuous on $(\VV, \|\cdot\| )$ almost surely.
\end{theorem}

When the kernel $u$ is symmetric, $\{\psi(\nu),\nu\in\VV\}$ is a
second-order Gaussian chaos and it is well known that we can take
%
\begin{equation}
\| \mu-\nu\| =\bigl(E\bigl(\psi(\mu)-\psi(\nu)\bigr)^{2}
\bigr)^{1/2}.\label{1.11}
\end{equation}

One of interests in studying permanental fields is to use them to
analyze families of continuous additive functionals.
We may think of a continuous additive functional of the Markov process
$X =
(\Om, \FF_{t}, X_t,\th_{t},P^x
)$ as
%
\begin{equation}
L^{\nu}_{t}:=\lim_{\ep\to0}\int
_S\int_0^t
\delta_{y,\ep}(X_s) \,ds \,d\nu(y),\label{0.2}
\end{equation}
where $\nu$ is a positive measure on $S$ and $\delta_{y,\ep}$ is an
approximate delta function at $y\in S$. More precisely, a family $A=\{
A_t; t\geq0\}$ of random
variables is called a continuous
additive functional
of
$X$
\label{CAF} if:
\begin{longlist}[(1)]
\item[(1)] $t\mapsto A_t$ is almost surely continuous and nondecreasing,
with
$A_0=0$ and $A_t=A_\ze$, for all $t\geq\ze$.
\item[(2)] $A_t$ is $\FF_t$ measurable.
\item[(3)]
$A_{t+s}=A_t+A_s\circ\th_t \mbox{ for all } s,t>0$ a.s.
\end{longlist}
(Details on the definition of $L^{\nu}_{t}$ are given in Section~\ref{sec-exist}.)

As in \cite{MR96} we relate permanental
fields and continuous additive functionals by a Dynkin type isomorphism
theorem. In Section~\ref{sec-it}, we obtain such a theorem relating $\{L^{\nu}_{\infty}\}$
and the
associated permanental field $\{\psi(\nu)\}$. Since the construction of
$\psi$ in Section~\ref{sec-exist} explores many properties of $\{
L^{\nu
}_{\infty}\}$, the further derivation of the isomorphism theorem is
relatively straightforward.

In Section~\ref{sec-it}, we introduce the measure
%
\begin{equation}
Q_{\phi}^{\rho}( F)= \int Q^{x,x}\bigl(F
L^{\phi}_{\infty} \bigr) \,d\rho(x), \label{nit4.4i}
\end{equation}
where $Q^{x,y}$ is given in (\ref{nit4.4}).

The next theorem is implied by Theorem~\ref{theo-ljrm}.

\begin{theorem}\label{theo-ljrmi} Let $X$ be a transient Borel right
process with
potential densities $u$, as described in Section~\ref{ss-caf}, and let
$\|\cdot\|$ be a proper norm for $u$. Let $\{\psi(\nu),\nu\in\RR
_{\|\cd
\|}\}$ be the associated $\al$-permanental field with kernel $u$. Then
for any $ \phi, \rho\in\mathcal{R}^{+}_{\|\cdot\|}$\vspace*{1pt} and all measures
$\{\nu_{j}\}\in\mathcal{R}_{\|\cdot\|}$, and all bounded
measurable functions
$F$ on $R^\infty$,
%
\begin{equation}
E Q_{\phi}^{\rho}\bigl( F \bigl( \psi(\nu_{i})+L_{\infty}^{\nu_{i}}
\bigr) \bigr)={1
\over\al
}E \bigl(\th^{\rho,\phi} F\bigl(\psi(
\nu_{i})\bigr)\bigr),\label{6.1i}
\end{equation}
where $\th^{\rho,\phi} $ is a random variable that has all moments finite.

[Here, we use the notation $F( f(x_{ i})):=F( f(x_{ 1}),f(x_{
2}),\ldots
)$, and the expectation of the $\{L_{\infty}^{\nu_{i}}\}$ are with respect
to $Q_{\phi}^{\rho}$, and of the $ \{\psi(\nu_{i})\}$ and $\{\th
^{\rho
,\phi} \}$ are with respect to $E$.]
\end{theorem}

It is easy to show that this isomorphism theorem implies that the
continuity of $\{\psi(\nu),\nu\in\VV\}$ on $(\VV, \|\cdot\|)$, implies
the continuity of $\{ L_{\infty}^{\nu},\nu\in\VV\}$ on $(\VV, \|
\cdot\|))$.
Extending this to the joint continuity of
$\{ L_{t}^{\nu},(\nu,t)\in\VV\times R^{+}\}$ on $( \VV\times
R^{+}, \|
\cdot\|\times|\cdot|)$ is considerably more difficult. We do this in
Section~\ref{sec-caf}.

Additional hypotheses are required to prove joint continuity
of $\{ L_{t}^{\nu},(\nu,t)\in\VV\times R^{+}\}$ in the most general
setting. However, these are satisfied by a simple sufficient condition
when the Markov process is a transient L\'evy processes.
Let $S=R^{d}$ and $ X$ be a L\'evy process killed at the end of an
independent exponential time, with characteristic function
%
\begin{equation}
Ee^{i\la X_{t}}=e^{-t\ka(\la)}\label{1.14}
\end{equation}
and potential density $u(x,y)=u(y-x)$.
We refer to $\ka$ as the characteristic exponent of $ X$.

We assume that
%
\begin{equation}
\|u\|_{2}<\infty\quad\mbox{and}\quad e^{- \operatorname{Re} \ka(\xi)}\qquad\mbox{is integrable on
$R^{d}$.} \label{expint}
\end{equation}
We say that $u$ is radially regular at infinity if
%
\begin{equation}
\frac{1}{ \tau(|\xi|) }\le\bigl|\hat u(\xi)\bigr|\le\frac{C }{ \tau(|\xi
|) }, \label{7.8q}
\end{equation}
where $ \tau(|\xi|)$ is regularly varying at infinity.
Note that
%
\begin{equation}
\hat u(\xi)={1 \over\ka(\xi)}. \label{1.19cor}
\end{equation}
For a measure $\nu$ on $R^{d}$, we define the measure $\nu_{h}$ by
%
\begin{equation}
\nu_{h}(A)=\nu(A-h).\label{1.19}
\end{equation}

\begin{theorem}\label{theo-1.1} Let
$X=\{X(t),t\in R^+\}$ be a L\'evy process in $R^d$ that is killed at
the end of an independent exponential time, with potential density
$u(x,y)=u(y-x)$. Assume that (\ref{expint}) holds and $\hat u$ is
radially regular. Let $\nu\in\mathcal{R}^{+}(X)$ and $\ga=|\hat
u|\ast
|\hat u|$. If
%
\begin{equation}
\int_{1}^{\infty} {(\int_{|\xi|\geq x}|\hat{\nu}(\xi)|^{2} \ga
(\xi) \,d\xi)
^{1/2}\over
x} \,dx<\infty,
\label{b101}
\end{equation}
then
$\{L_{t}^{\nu_x}, (x,t)\in R^d\times R_{+}\}$ is continuous $P^{y}$
almost surely for all $y\in R^d$.

In addition
%
\begin{equation}
\limsup_{\delta\to0}\mathop{\sup_{|x-y|\le\delta}}_{ x,y\in[0,1]^d} \frac{L^{\nu_{x}}_{t}-L^{\nu_{y}}_{t}}{\omega(\delta) }
\le C\qquad \mbox{a.s.},\label{15.1ai}
\end{equation}
where
%
\begin{equation}
\omega(\delta)= \varphi(\delta)\log1/\delta+\int_{0}^{\delta}
\frac{ \varphi(u)}{u} \,du\label{1.22}
\end{equation}
and
%
\begin{equation}\qquad\hspace*{6pt}
\varphi(\delta) = \biggl( |\delta|^{2} \int_{|\xi|\le
1/|\delta|}|
\xi |^{2}\bigl|\hat\nu(\xi)\bigr|^{2} \ga(\xi) \,d\xi+\int
_{|\xi|\ge1/|\delta
|} \bigl|\hat\nu (\xi)\bigr|^{2} \ga(\xi) \,d\xi
\biggr)^{1/2}.\label{1.23}
\end{equation}
\end{theorem}

\begin{example}\label{ex-1.1}
1. If $\tau(|\xi|)$ is regularly varying at infinity with index
greater than $ d/2 $ and less than $d$ and
%
\begin{equation}
\bigl|\hat\nu(\xi)\bigr| \le C \frac{\tau(|\xi|)}{ |\xi|^{d} (\log|\xi|)^{
3/2+\ep}} \qquad\mbox{as } |\xi| \rar
\infty\label{amp2}
\end{equation}
for some constant $C>0$ and any $\ep>0$, then $\{L_{t}^{\nu_x},
(x,t)\in R^d\times R_{+}\}$ is continuous~$P^{x}$ almost surely.
\begin{longlist}[2.]
\item[2.] If
%
\begin{equation}
\tau\bigl(|\xi|\bigr)=\frac{|\xi|^{2}} {(\log|\xi|)^{a}}\qquad \mbox{for $a\ge 0$ and all $|\xi|$ sufficiently
large}\label{1.25}
\end{equation}
and
%
\begin{equation}
\bigl|\hat\nu(\xi)\bigr| \le C \frac{\tau(|\xi|)}{ |\xi|^{2} (\log|\xi|)^{
2+\ep}} \qquad\mbox{as } |\xi| \rar
\infty\label{amp2q}
\end{equation}
for some constant $C>0$ and any $\ep>0$, then $\{L_{t}^{\nu_x},
(x,t)\in R^2\times R_{+}\}$ is continuous $P^{x}$ almost surely. This
extends the result for Brownian motion in $R^{2}$ (in which case $a=0$)
that is given in \cite{MR96}, Theorem~1.6.

\item[3.]
If $\tau(|\xi|)$ is regularly varying at infinity with index $d/2<\al
<d$ and
%
\begin{equation}
\bigl|\hat\nu(\xi)\bigr| \le\frac{1}{\vartheta(|\xi|)},
\end{equation}
where $\vartheta(|\xi|)$ is regularly varying at infinity with index
$\bb$ and $\al+\bb>d $, then
there exists a constant $C>0$, such that for almost every $t $,
%
\begin{equation}
\limsup_{\delta\to0}\mathop{\sup_{|x-y|\le\delta}}_{ x,y\in[0,1]^d} \frac{L^{\nu_{x}}_{t}-L^{\nu_{y}}_{t}}{\varrho(\delta)\log
1/\delta}
\le C \qquad\mbox{a.s.},\label{15.1aiz}
\end{equation}
where
%
\begin{equation}
\varrho(\delta)\sim C \bigl(\delta^{-d }\tau(1/\delta)\vartheta (1/
\delta)\bigr)^{-1} \qquad\mbox {as $\delta\to0$},
\end{equation}
is regularly varying at zero with index $\al+\bb-d$.

\item[4.]
If $d=2$ and $\tau(|\xi|)$ is as given in (\ref{1.25}) and
%
\begin{equation}
\bigl|\hat\nu(\xi)\bigr| \le\frac{1}{\vartheta(|\xi|)},
\end{equation}
where $\vartheta(|\xi|)$ is regularly varying at infinity with index
$\bb>0$, then there exists a constant $C>0$, such that for almost every
$t $
%
\begin{equation}
\limsup_{\delta\to0}\mathop{\sup_{|x-y|\le\delta}}_{ x,y\in[0,1]^d} \frac{L^{\nu_{x}}_{t}-L^{\nu_{y}}_{t}}{\varrho(\delta)\log
1/\delta}
\le C \qquad\mbox{a.s.},\label{15.1aizc}
\end{equation}
where
%
\begin{equation}
\varrho(\delta)\sim C \bigl(\delta^{-2}\tau(1/\delta)\vartheta (1/
\delta)\bigr)^{-1} (\log 1/\delta)^{1/2} \qquad\mbox{as $\delta
\to0$},
\end{equation}
is regularly varying at zero with index $ \bb$.
\end{longlist}
\end{example}

Continuous additive functionals of L\'evy processes are studied in
Section~\ref{sec-levy}.

\section{Markov loops and the existence of permanental fields}\label
{sec-exist}

So far, a permanental field is defined as a process with a certain
moment structure. In this section, we show that a permanental field
with kernel $u$ can be realized in terms of continuous additive
functionals of a Markov process $X$ with potential density $u$.

\subsection{Continuous additive functionals}\label{ss-caf}

Let $S$ a be locally compact set with a countable base.
Let $X =
(\Om, \FF_{t}, X_t,\th_{t},P^x
)$ be a transient Borel right process with state space $S$, and jointly
measurable transition densities $p_{t}(x,y)$ with respect to some $\si
$-finite measure $m$ on $S$.
We assume that the potential densities
%
\begin{equation}
u(x,y)=\int_{0}^{\infty}p_{t}(x,y) \,dt
\end{equation}
are finite off the diagonal, but allow them to be infinite on the
diagonal. We also assume that $\sup_{x}\int_{\delta}^{\infty
}p_{t}(x,x) \,dt<\infty
$ for each $\delta>0$. We do not require that $p_{t}(x,y)$ is symmetric.

We assume furthermore that $0<p_{t}(x,y)<\infty$ for all $0<t<\infty
$ and $
x,y\in S$, and that there exists another right process $\widehat X$ in
duality with $X$, relative to the measure $m$, so that its transition
probabilities
$\widehat P_{t}(x,dy)=p_{t}(x,y) m(dy)$. These conditions allow us to use
material on bridge measures in \cite{FPY} in the construction of the
loop measure in Section~\ref{ss-lm}.

Let $\{A_t,t\in R^{+}\}$ be a positive continuous additive functional
of $X$. The 0-potential of $\{A_t,t\in R^{+}\}$ is defined to be
%
\begin{equation}
u^{0}_{A}(x)=E^{x}(A_{\infty}).\label{potz.1}
\end{equation}
If $\{A_t,t\in R^{+}\}$ and $\{B_{t}, t\in R^{+}\}$ are two continuous
additive functionals of $X$, with $u^{0}_{A}=u^{0}_{B}<\infty$, then
$\{
A_t, t\in R^{+}\}=\{ B_{t}, t\in R^{+}\}$ a.s.
(see, e.g., \cite{S}, Theorem~36.10). This can also be seen directly by
noting that the properties of a continuous additive functional given in
its definition and the
Markov property imply that $M_{t}=A_t- B_{t}$ is a continuous
martingale of bounded variation, and consequently is a constant,
\cite{RY}, Chapter~IV, Proposition~1.2, which is zero in this case because $M_{0}=0$.

When $\{A_t,t\in R^{+}\}$ is a positive continuous additive functional
with $0$-potential $u^{0}_{A}$ that is the potential of a $\sigma
$-finite measure $\nu$, that is when
%
\begin{equation}
E^{x}(A_{\infty})=\int u(x,y) \,d\nu(y),\label{potz.2}
\end{equation}
we write $A_{t}=L^\nu_t$ and refer to $\nu$ as the Revuz measure of $A_t$.

It follows from
\cite{Revuz}, Section~V.6, that a $\sigma$-finite measure is the Revuz measure
of a continuous additive functional of a Markov process $X$ with
potential density $u$ if and
only if
%
\begin{equation}
U\nu(x):=\int u(x,y) \,d\nu(y)<\infty\qquad\mbox{for each $x\in S$}
\end{equation}
and $\nu$ does not charge any semi-polar set.
We denote by $\mathcal{R}^{+}(X)$, or $\mathcal{R}^{+}$ when $X$ is
understood, the set of positive bounded Revuz measures. We use
$\mathcal
{R}$ for the set of measures of the form $\nu=\nu_{1}-\nu_{2}$ with $
\nu_{1},\nu_{2}\in\mathcal{R}^{+}$, and we write $L^\nu_{t}=L^{
\nu
_{1}}_{t}-L^{ \nu_{2}}_{t}$. The comments above show that this is well
defined. Throughout this paper, we only consider measures in $\mathcal{R}$.

\subsection{Loop measure}\label{ss-lm}

It follows from the assumptions in the first two paragraphs of Section~\ref{ss-caf} that, as in \cite{FPY}, for all $0<t<\infty$ and
$x,y\in S$,
there exists a finite measure $Q_{t}^{x,y}$ on $\mathcal{F}_{t^{-}}$,
of total mass $p_{t}(x,y)$, such that
%
\begin{equation}
Q_{t}^{x,y}(1_{\{\ze>s\}} F_{s})=P^{x}
\bigl( F_{s} p_{t-s}(X_{s},y)
\bigr)\label{10.1}
\end{equation}
for all $F_{s} \in\mathcal{F}_{s}$ with $s<t$. (In this paper, we use
the letter $Q$ for measures which are not necessarily of mass $1$, and
reserve the letter $P$ for probability measures.)

We use the canonical representation of $X$ in which $\Om$ is the set
of right continuous paths $\omega$ in $S_{\De}=S\cup\De$ with $\De
\notin
S$, and is such that $\omega(t)=\De$ for all $t\geq\ze=\inf\{t>0
|\omega
(t)=\De\}$. Set $X_{t}(\omega)=\omega(t)$. We define a $\si
$-finite measure
$\mu$
on $(\Om, \mathcal{F})$ by
%
\begin{equation}
\mu(F)=\int_{0}^{\infty}{1\over t}\int
Q_{t}^{x,x}(F\circ k_{t}) \,dm (x) \,dt\label{ls.3}
\end{equation}
for all $\mathcal{F}$ measurable functions $F$ on $\Om$.
Here, $k_{t}$ is the killing operator defined by $k_{t}\omega
(s)=\omega(s)$
if $s<t$ and $k_{t}\omega(s)=\De$ if $s\geq t$, so that
$k_{t}^{-1}\mathcal
{F}\subset\mathcal{F}_{t^{-}}$.
We call $\mu$ the loop measure of $X$, because when $X$ has continuous
paths, $\mu$ is concentrated on the set of continuous
loops with a distinguished starting point (since $Q_{t}^{x,x}$ is
carried by loops starting at $x$). It can be shown that $\mu$ is
invariant under ``loop rotation,'' and $\mu$ is often restricted to the
``loop rotation'' invariant sets. We do not pursue these ideas in this paper.

As usual, if $F$ is a function, we often write $\mu(F)$ for $\int F
\,d\mu$. [We already used this notation in (\ref{10.1}).]

We explore some properties of the loop measure $\mu$. [Recall
the definition of $ \mathcal{R}_{\|\cdot\|}$ in the paragraph
containing (\ref{rev1i}).]

\begin{lemma}\label{lem-ls} Let $k\geq2$, and assume that $\nu
_{j}\in
\mathcal{R}_{\|\cdot\|}$ for all $j=1,\ldots,k$. Then
%
\begin{eqnarray}\label{ls.4}
&&\mu\Biggl( \prod_{j=1}^{k}L^{\nu_{j}}_{\infty}
\Biggr)
\nonumber
\\[-8pt]
\\[-8pt]
\nonumber
&&\qquad ={1 \over k}\sum_{\pi\in\mathcal{P}_{k}}\int
u(y_{1},y_{2})\cdots u(y_{k-1},y_{k})u(y_{k},y_{1})
\prod_{j=1}^{k} \,d\nu_{\pi
(j)}(y_{j}),
\end{eqnarray}
where $\mathcal{P}_{k}$ denotes the set of permutations of $[1,k]$.
Equivalently,
%
\begin{eqnarray}\label{8.2t}
&& \mu\Biggl( \prod_{j=1}^{k}L^{\nu_{j}}_{\infty}
\Biggr)\nonumber\\
&&\qquad=\sum_{\pi\in\mathcal
{P}_{k-1}}\int \Biggl(\int u(x,
y_{1})u(y_{1},y_{2})
\\
&&\hspace*{69pt}\qquad \cdots u(y_{k-2},y_{k-1})u(y_{k-1},x) \prod
_{j=1}^{k-1} \,d\nu_{\pi(j)}(y_{j})
\Biggr) \,d\nu_{k}(x).\nonumber
\end{eqnarray}
\end{lemma}

When $k=1$, the formula in (\ref{ls.4}) gives
%
\begin{equation}
\mu\bigl( L^{\nu}_{\infty}\bigr)=\int u(y,y) \,d\nu(y).
\end{equation}
Obviously, this is infinite when $u(y,y)=\infty$.

\begin{pf*}{Proof of Lemma \ref{lem-ls}} We first assume that all the $\nu_{j}$ are positive
measures. Note that for all $j=1,\ldots,k$
%
\begin{equation}
L^{\nu_{j}}_{\infty} \circ k_{t}=L^{\nu_{j}}_{t}.\label{10.5}
\end{equation}
Therefore,
%
\begin{eqnarray}\label{ls.1z}
Q_{t}^{x,x}\Biggl( \Biggl(\prod
_{j=1}^{k}L^{\nu_{j}}_{\infty}\Biggr)
\circ k_{t}\Biggr)&=&Q_{t}^{x,x}\Biggl( \prod
_{j=1}^{k}L^{\nu_{j}}_{t}
\Biggr)
\nonumber
\\
&=& Q_{t}^{x,x}\Biggl( \prod_{j=1}^{k}
\int_{0}^{t} \,dL^{\nu
_{j}}_{r_{j}}
\Biggr)
\\
&=& \sum_{\pi\in\mathcal{P}_{k}} Q_{t}^{x,x}
\biggl( \int_{0\leq
r_{1}\leq
\cdots\leq r_{k}\leq t} \,dL^{\nu_{\pi(1)}}_{r_{1}}\cdots
dL^{\nu
_{\pi
(k)}}_{r_{k}}\biggr).\nonumber
\end{eqnarray}
We use the following technical lemma.

\begin{lemma}\label{lem-finq} Let $\nu_{j}\in\mathcal{R}^{+}$ for all
$j=1,\ldots
,k$. Then for all $t\in R^{+}$
%
\begin{eqnarray}\label{ls.1}
&& Q_{t}^{x,y}\biggl( \int_{0\leq r_{1}\leq\cdots\leq r_{k}\leq t}
\,dL^{\nu
_{1}}_{r_{1}}\cdots dL^{\nu_{k}}_{r_{k}}\biggr)
\nonumber\\
&&\qquad= \int_{0\leq r_{1}\leq\cdots\leq r_{k}\leq t}\int p_{r_{1}}(x,y_{1})
p_{r_{2}-r_{1}}(y_{1},y_{2})
\\
&&\hspace*{73pt}\qquad\quad\cdots p_{r_{k}-r_{k-1}}(y_{k-1},y_{k})
p_{t-r_{k}}(y_{k},y) \prod_{j=1}^{k}
\,d\nu_{j}(y_{j}) \,dr_{j}.
\nonumber
\end{eqnarray}
\end{lemma}

\begin{pf} We prove this by induction on $k$. The case $k=1$
follows from \cite{FPY}, Lemma~1. Assume we have proved (\ref{ls.1})
for all $1\le j\le k-1$.
We write
%
\begin{eqnarray}\label{ls.1jj}
&& Q_{t}^{x,y}\biggl( \int_{0\leq r_{1}\leq\cdots\leq r_{k}\leq t}
\,dL^{\nu
_{1}}_{r_{1}}\cdots dL^{\nu_{k}}_{r_{k}}\biggr)
\nonumber
\\[-8pt]
\\[-8pt]
\nonumber
&&\qquad=Q_{t}^{x,y}\biggl( \int_{0}^{t}H_{r_{k}}
\,dL^{\nu_{k}}_{r_{k}}\biggr),
\end{eqnarray}
where
%
\begin{equation}
H_{r_{k}}= \int_{0\leq r_{1}\leq\cdots\leq r_{k}} \,dL^{\nu
_{1}}_{r_{1}}
\cdots dL^{\nu_{k-1}}_{r_{k-1}}.
\end{equation}
Clearly, $H_{r_{k}}$ is continuous in $r_{k}$.
It follows from \cite{FPY}, Proposition~3, that
%
\begin{eqnarray}\label{nl1}
&& Q_{t}^{x,y}\biggl( \int_{0}^{t}H_{r_{k}}
\,dL^{\nu_{k}}_{r_{k}}\biggr)
\nonumber
\\[-8pt]
\\[-8pt]
\nonumber
&&\qquad=Q_{t}^{x,y}\biggl(
\int_{0}^{t}{Q_{r_{k}}^{x,X_{r_{k}}}( H_{r_{k}}) \over
p_{r_{k}}(x,X_{r_{k}})}
\,dL^{\nu_{k}}_{r_{k}}\biggr).
\end{eqnarray}
Using \cite{FPY}, Lemma~1, again, we see that
%
\begin{eqnarray} \label{nl2}
&&Q_{t}^{x,y}\biggl( \int_{0}^{t}
{Q_{r_{k}}^{x,X_{r_{k}}}( H_{r_{k}}) \over
p_{r_{k}}(x,X_{r_{k}})} \,dL^{\nu_{k}}_{r_{k}}\biggr)
\nonumber\\
&&\qquad =\int_{0}^{t} \biggl(\int
p_{r_{k}}(x,y_{k})p_{t-r_{k}}(y_{k},y)
{Q_{r_{k}}^{x,y_{k}}( H_{r_{k}}) \over p_{r_{k}}(x,y_{k})} \,d\nu _{k}(y_{k})\biggr)
\,dr_{k}
\\
&&\qquad =\int_{0}^{t} \biggl(\int
p_{t-r_{k}}(y_{k},y) Q_{r_{k}}^{x,y_{k}}(
H_{r_{k}}) \,d\nu_{k}(y_{k})\biggr)
\,dr_{k}.
\nonumber
\end{eqnarray}
Using (\ref{ls.1jj}) and (\ref{ls.1}) for $k-1$, we see that it holds
for all $1\le j\le k-1$.
\end{pf}

\begin{pf*}{Proof of Lemma~\ref{lem-ls} continued} Combining
(\ref
{ls.1}) with (\ref{ls.1z}), we obtain
%
\begin{eqnarray}\label{ls.2y}
&&Q_{t}^{x,x}\Biggl( \prod
_{j=1}^{k}L^{\nu_{j}}_{\infty}\circ
k_{t}\Biggr) \nonumber
\\
&&\qquad=\sum_{\pi\in\mathcal{P}_{k}} \int_{0\leq r_{1}\leq\cdots\leq
r_{k}\leq t} \int
p_{r_{1}}(x,y_{1}) p_{r_{2}-r_{1}}(y_{1},y_{2})
\nonumber
\\[-8pt]
\\[-8pt]
\nonumber
&&\hspace*{94pt}\qquad\quad\cdots p_{r_{k}-r_{k-1}}(y_{k-1},y_{k})
p_{t-r_{k}}(y_{k},x)\\
&&\hspace*{94pt}\qquad\quad{}\times \prod_{j=1}^{k}
\,d\nu_{\pi(j)}(y_{j}) \,dr_{j}.
\nonumber
\end{eqnarray}
Therefore,
%
\begin{eqnarray}\label{ls.2}
&& \int Q_{t}^{x,x}\Biggl( \prod
_{j=1}^{k}L^{\nu_{j}}_{\infty}\circ
k_{t}\Biggr) \,dm (x)
\nonumber\\
&&\qquad=\sum_{\pi\in\mathcal{P}_{k}} \int_{0\leq r_{1}\leq\cdots\leq
r_{k}\leq t} \int
p_{r_{2}-r_{1}}(y_{1},y_{2})
\nonumber
\\[-8pt]
\\[-8pt]
\nonumber
&&\hspace*{95pt}\qquad\quad\cdots p_{r_{k}-r_{k-1}}(y_{k-1},y_{k})
p_{r_{1}+t-r_{k}}(y_{k},y_{1})\\
&&\hspace*{103pt}\qquad{}\times \prod
_{j=1}^{k} \,d\nu_{\pi(j)}(y_{j})
\,dr_{j},
\nonumber
\end{eqnarray}
since
%
\begin{equation}
\int p_{r_{1}}(x,y_{1})p_{t-r_{k}}(y_{k},x)
\,dm(x)=p_{r_{1}+t-r_{k}}(y_{k},y_{1}).
\end{equation}
It follows from (\ref{ls.3}) and (\ref{ls.2})
that
%
\begin{eqnarray}\label{ls.2m}
&&\mu\Biggl( \prod_{j=1}^{k}L^{\nu_{j}}_{\infty}
\Biggr)
\nonumber
\\
&&\qquad=\sum_{\pi\in\mathcal{P}_{k}} \int_{0}^{\infty}
{1\over t}\Biggl(\int_{0\leq
r_{1}\leq\cdots\leq r_{k}\leq t} \int
p_{r_{2}-r_{1}}(y_{1},y_{2})
\nonumber
\\[-8pt]
\\[-8pt]
\nonumber
&&\hspace*{130pt}\qquad\quad \cdots p_{r_{k}-r_{k-1}}(y_{k-1},y_{k})
p_{r_{1}+t-r_{k}}(y_{k},y_{1})\\
&&\hspace*{210pt}\qquad{}\times \prod
_{j=1}^{k} \,d\nu_{\pi(j)}(y_{j})
\,dr_{j}\Biggr) \,dt.\nonumber
\end{eqnarray}

We make the change of variables $(t,r_{2},\ldots,r_{k})\to
s_{1}=r_{1}+t-r_{k}, s_{2}=r_{2}-r_{1},\ldots, s_{k}=r_{k}-r_{k-1}$,
and integrate on $r_{1}$ to obtain\vspace*{1pt}
%
\begin{eqnarray}\label{ls.5c}
\hspace*{-4pt}&&\mu\Biggl( \prod_{j=1}^{k}L^{\nu_{j}}_{\infty}
\Biggr)
\nonumber
\\[1pt]
\hspace*{-6pt}&&\qquad=\sum_{\pi\in\mathcal{P}_{k}}\int{1\over s_{1}+\cdots+s_{k}}\Biggl(
\int p_{s_{2}}(y_{1},y_{2})
\nonumber
\\[1pt]
\hspace*{-6pt}&&\hspace*{111pt}\qquad\quad\cdots p_{s_{k}}(y_{k-1},y_{k})
p_{s_{1}}(y_{k},y_{1})
\nonumber\\
\hspace*{-6pt}&&\hspace*{186pt}{}\times  \prod
_{j=1}^{k} \,d\nu_{\pi(j)}(y_{j})
\Biggr) \biggl(\int_{0}^{s_{1}} 1 \,dr_{1}
\biggr)\prod_{j=1}^{k} \,ds_{j}
\\[1pt]
\hspace*{-6pt}&&\qquad=\sum_{\pi\in\mathcal{P}_{k}}\int{s_{1} \over s_{1}+\cdots+s_{k}} \int
p_{s_{2}}(y_{1},y_{2})
\nonumber
\\
\hspace*{-6pt}&&\hspace*{105pt}\qquad\quad\cdots p_{s_{k}}(y_{k-1},y_{k})
p_{s_{1}}(y_{k},y_{1}) \prod
_{j=1}^{k} \,d\nu_{\pi(j)}(y_{j})
\,ds_{j}.\nonumber
\end{eqnarray}
Set
%
\begin{eqnarray}\label{spe}
&& f(s_{1}, s_{2}, \ldots, s_{k})
\nonumber
\\[-8pt]
\\[-8pt]
\nonumber
&&\qquad=\sum_{\pi\in\mathcal{P}_{k}}\int p_{s_{2}}(y_{1},y_{2})
\cdots p_{s_{k}}(y_{k-1},y_{k}) p_{s_{1}}(y_{k},y_{1})
\prod_{j=1}^{k} \,d\nu _{\pi(j)}(y_{j}),
\end{eqnarray}
and note that because of the sum over all permutations
%
\begin{equation}
f(s_{1}, s_{2}, \ldots, s_{k})=f(s_{2},
s_{3}, \ldots, s_{1}).\label{spf}
\end{equation}
Using (\ref{spf}) after a simple change of variables, we see from
(\ref
{ls.5c}) that
%
\begin{eqnarray}\label{8.2q}
\mu\Biggl( \prod_{j=1}^{k}L^{\nu_{j}}_{\infty}
\Biggr)&=& \int{s_{1} \over
s_{1}+\cdots+s_{k}}f(s_{1}, s_{2},
\ldots, s_{k}) \prod_{j=1}^{k}
\,ds_{j}
\nonumber\\
&=& \int{s_{2} \over s_{2}+s_{3}+\cdots+s_{1}}f(s_{2}, s_{3}, \ldots,
s_{1}) \prod_{j=1}^{k}
\,ds_{j}
\\
&=& \int{s_{2} \over s_{1}+s_{2}+\cdots+s_{k}}f(s_{1}, s_{2}, \ldots,
s_{k}) \prod_{j=1}^{k}
\,ds_{j}.
\nonumber
\end{eqnarray}
Similarly, we see that for all $1\leq j\leq k$
%
\begin{equation}
\mu\Biggl( \prod_{j=1}^{k}L^{\nu_{j}}_{\infty}
\Biggr)=\int{s_{j} \over
s_{1}+s_{2}+\cdots+s_{k}}f(s_{1}, s_{2}, \ldots,
s_{k}) \prod_{j=1}^{k}
\,ds_{j}. \label{}
\end{equation}
Therefore,
%
\begin{eqnarray}\label{8.2}
\mu\Biggl( \prod_{j=1}^{k}L^{\nu_{j}}_{\infty}
\Biggr)&=&{1 \over k} \int {s_{1}+\cdots+s_{k} \over s_{1}+\cdots+s_{k}}f(s_{1},
s_{2}, \ldots, s_{k}) \prod_{j=1}^{k}
\,ds_{j}
\nonumber
\\
&=&{1 \over k} \int f(s_{1}, s_{2}, \ldots,
s_{k}) \prod_{j=1}^{k}
\,ds_{j}
\nonumber
\\
&=&{1 \over k}\sum_{\pi\in\mathcal{P}_{k}}\int\int
p_{s_{2}}(y_{1},y_{2})\cdots p_{s_{k}}(y_{k-1},y_{k})p_{s_{1}}(y_{k},y_{1})
\\
&&\hspace*{53pt}{}\times \prod_{j=1}^{k} \,d\nu_{\pi(j)}(y_{j})
\,ds_{j}\nonumber
\\
&=&{1 \over k}\sum_{\pi\in\mathcal{P}_{k}}\int
u(y_{1},y_{2})\cdots u(y_{k-1},y_{k})u(y_{k},y_{1})
\prod_{j=1}^{k} \,d\nu_{\pi(j)}(y_{j}).\nonumber
\end{eqnarray}
It follows from the hypothesis that $\|\cdot\|$ is a proper norm that
the integrals in (\ref{8.2}) are finite; consequently, the equalities
in (\ref{8.2}) hold for all $\nu\in\RR_{\|\cdot\|}$ (i.e., measures
that are not necessarily positive). This is (\ref{ls.4}).

To obtain (\ref{8.2t}), we note that because
we are permuting $k$ points on a circle, for $k\ge2$, we can write
(\ref{ls.4}) as
%
\begin{eqnarray}\label{8.2tq}
&&\mu\Biggl( \prod_{j=1}^{k}L^{\nu_{j}}_{\infty}
\Biggr)\nonumber\\
&&\qquad=\sum_{\pi\in\mathcal
{P}_{k-1}}\int \Biggl(\int u(x,
y_{1})u(y_{1},y_{2})
\\
&&\qquad\hspace*{68pt}{} \cdots u(y_{k-2},y_{k-1})u(y_{k-1},x) \prod
_{j=1}^{k-1} \,d\nu_{\pi(j)}(y_{j})
\Biggr) \,d\nu_{k}(x).
\nonumber
\end{eqnarray}
\upqed\end{pf*}
\noqed\end{pf*}

\begin{remark}\label{rem-3.1} Note that in the course of the proof
of Lemma~\ref{lem-ls} we show that (\ref{ls.4}) and (\ref{8.2t}) hold
for all measures in $\RR^{+}$.
\end{remark}

For use in the next section, note that when $k=1$, (\ref{ls.2}) takes
the form
%
\begin{equation}
\int Q_{t}^{x,x}\bigl( L^{\nu}_{\infty}\circ
k_{t}\bigr) \,dm (x)=t \int p_{t}(y,y ) \,d\nu(y ).
\label{ls.2mq}
\end{equation}
Using the fact that $1_{\{\ze>\delta\}}\circ k_{t}=1$ if $t>\delta
$, and $0$
if $t\leq\delta$, we see that
%
\begin{equation}
\mu\bigl( 1_{\{\ze>\delta\}}L_{\infty}^{\nu}\bigr)=\int
_{\delta
}^{\infty}\int p_{t}(y,y ) \,d\nu(y )
\,dt,\label{ls.2mm}
\end{equation}
which is finite by our assumptions that $\sup_{x}\int_{\delta
}^{\infty
}p_{t}(x,x) \,dt<\infty$ for each $\delta>0$ and $\nu$ is a finite measure.

\subsection{Loop soup}
\label{ss-ls}

Let $\mathcal{L}_{\al}$ be the Poisson point process on $\Om$ with
intensity measure $\al\mu$. Note that $\mathcal{L}_{\al}$ is a random
variable; each realization of $\mathcal{L}_{\al}$ is a countable subset
of $\Om$. To be more specific, let
%
\begin{equation}
N(A):=\# \{\mathcal{L}_{\al}\cap A \},\qquad A\subseteq\Om.\label{po.1}
\end{equation}
Then for any disjoint measurable subsets $A_{1},\ldots,A_{n}$ of $\Om
$, the random variables $N(A_{1}),\ldots,N(A_{n})$, are independent
and $N(A)$ is a Poisson random variable with parameter $\al\mu(A)$, that is,
%
\begin{equation}
P_{\mathcal{L}_{\al}}\bigl(N(A)=k\bigr)={(\al\mu(A))^{k} \over k!}e^{-\al\mu
(A)}.\label{po.2}
\end{equation}

The Poisson point process $\mathcal{L}_{\al}$ is called the loop soup
of the Markov process $X$. For $\nu\in\mathcal{R}_{\|\cdot\|}$, we define
%
\begin{equation}
\widetilde \psi(\nu)=\lim_{\delta\rar0}\widehat L_{\delta}^{\nu
},
\label{ls.10a}
\end{equation}
where
%
\begin{equation}
\widehat L_{\delta}^{\nu} =\biggl(\sum
_{\omega\in\mathcal{L}_{\al}} 1_{\{\ze(\omega)>\delta\}
} L_{\infty}^{\nu
}(
\omega)\biggr)-\al\mu\bigl( 1_{\{\ze>\delta\}}L_{\infty}^{\nu}\bigr).
\label{ls.10}
\end{equation}
As noted following (\ref{ls.2mm}), $\mu( 1_{\{\ze>\delta\}
}L_{\infty}^{\nu})$
is finite for all $\delta>0$.
We show in Theorem~\ref{theo-9.1} that the limit (\ref{ls.10a})
converges in all $L^{p}$, even though each term in (\ref{ls.10}) has an
infinite limit as $\delta\rar0$.

The terms loop soup and ``loop soup local time'' are used in \cite{LF,LW}, and \cite{LL}, Chapter~9. In \cite{LeJan}, they are referred to,
less colorfully albeit more descriptively, as
Poisson ensembles of Markov loops, and occupation fields of Poisson
ensembles of Markov loops.

The next theorem contains Theorem~\ref{theo-9.1i}. It is
given for symmetric kernels in~\cite{LeJan1}, Theorem~9. (In which
case, when $\al=1/2$, the permanental process is a second-order
Gaussian chaos.)

\begin{theorem}\label{theo-9.1} Let $X$ be a transient Borel right
process with state space $S$ and potential density $u(x,y)$, $x,y\in
S$, as described in the beginning of this section. Then for $\nu\in
\mathcal{R}_{\|\cdot\|}$, the limit (\ref{ls.10a}) converges in all
$L^{p}$ and $\{\widetilde \psi(\nu),\nu\in\mathcal{R}_{\|\cdot\|
}\}$ is an
$\al$-permanental field with kernel $u(x,y)$.
\end{theorem}

\begin{pf}
By the master formula for Poisson processes \cite{K}, (3.6),
%
\begin{eqnarray}\label{ls.11}
&& E_{\mathcal{L}_{\al}}\bigl(e^{ \sum_{j=1}^{n}z_{j}\widehat L_{\delta
_{j}}^{\nu_{j}}
}\bigr)
\nonumber
\\[-1pt]
\\[-15pt]
\nonumber
&&\qquad=\exp\Biggl(\al\Biggl(\int_{\Om}\Biggl( e^{\sum_{j=1}^{n}z_{j}1_{\{\ze>\delta
_{j}\}}
L^{\nu_{j}}_{\infty} } -
\sum_{j=1}^{n}z_{j}1_{\{\ze>\delta_{j}\}}
L^{\nu_{j}}_{\infty} -1\Biggr) \,d\mu(\omega) \Biggr)\Biggr).
\nonumber
\end{eqnarray}
Differentiating each side of (\ref{ls.11}) with respect to
$z_{1},\ldots
, z_{n}$ and then setting $z_{1},\ldots, z_{n}$ equal to zero, we see that
%
\begin{equation}
E_{\mathcal{L}_{\al}}\Biggl( \prod_{j=1}^{n}
\widehat L_{\delta_{j}}^{\nu
_{j}} \Biggr)= \sum
_{\bigcup_{i } B_{i}=[1,n], |B_{i}|\geq2} \prod_{i } \al\mu\biggl(
\prod_{j\in
B_{i}}1_{\{\ze>\delta_{j}\}}L^{\nu_{j}}_{\infty}
\biggr),\label{mom.11}
\end{equation}
where the sum is over all partitions $B_{1},\ldots, B_{n}$ of $[1,n]$
with all $|B_{i}|\geq2$.

The right-hand side of (\ref{mom.11}) can be written as a sum of terms
involving only positive measures, to which the monotone convergence
theorem can be applied. Using (\ref{8.2t}), we then see that the
right-hand side has a limit as the $\delta_{j}\rar0$ and this limit is
the same as the right-hand side of (\ref{a.3}). Applying this with
$\prod_{j=1}^{n}\widehat L_{\delta_{j}}^{\nu_{j}}$ replaced by
$(\widehat
L_{\delta
}^{\nu}-\widehat L_{\delta'}^{\nu})^{n}$, for arbitrary integer $n$,
shows that
the limit (\ref{ls.10a}) exists in all $L^{p}$.
\end{pf}

\begin{remark}
If we let $\alpha$ vary, we get a field-valued process with independent
stationary increments. This property is inherited from the analogous
property of the loop soup.
\end{remark}

\section{Isomorphism theorem}
\label{sec-it}

In this section, we obtain an isomorphism theorem that relates
permanental fields and continuous additive functionals. To begin, we
consider properties of several measures on the probability space of
$X$. Recall that $u$ denotes the 0-potential density of $X$.

Let $Q^{x,y}$ denote the $\sigma$-finite measure defined by
%
\begin{equation}
Q^{x,y}( 1_{\{\ze>s \}} F_s)= P^x\bigl(
F_s u(X_s,y )\bigr) \qquad\mbox{for all $F_s
\in b\mathcal{F}^{0}_s$}, \label{nit4.4}
\end{equation}
where $\mathcal{F}^{0}_s$ is the $\sigma$-algebra generated by $\{
X_{r}, 0\leq r\leq s\}$.

\begin{lemma}\label{lem-q} For all $x,y$
%
\begin{equation}
Q^{x,y} (F)=\int_{0}^{\infty}
Q_{t}^{x,y}(F\circ k_{t}) \,dt,\qquad  F\in b
\mathcal{F}^{0}.\label{ls.3n}
\end{equation}
\end{lemma}

\begin{pf}
To obtain (\ref{ls.3n}), it suffices to prove it for $F$ of the form
$\{1_{\{\ze>s \}} F_s\}$ for all $F_s\in b\mathcal{F}^{0}_s$.
Since $1_{\{\ze>s \}}\circ k_{t}=1_{\{s>t \}}1_{\{\ze>s \}} $,
%
\begin{eqnarray}\label{ls.3p}
&&\int_{0}^{\infty} Q_{t}^{x,y}
\bigl((1_{\{\ze>s \}} F_s)\circ k_{t}\bigr) \,dt\nonumber\\
&&\qquad= \int
_{s}^{\infty} Q_{t}^{x,y}(1_{\{\ze>s \}}
F_s) \,dt
\nonumber
\\[-8pt]
\\[-8pt]
\nonumber
&&\qquad=\int_{s}^{\infty} P^{x}\bigl(
F_s p_{t-s}(X_{s},y) \bigr) \,dt
\\
&&\qquad= P^x\bigl( F_s u(X_s,y )\bigr)=Q
^{x,y}(1_{\{\ze>s \}} F_s ), \nonumber
\end{eqnarray}
where the second and third equalities follow from (\ref{10.1}) and
interchanging the order of integration and the final equation by (\ref
{nit4.4}).
\end{pf}

We have the following formula for the moments of $\{L_{\infty}^{\nu
},\nu
\in\RR^{+}\}$ under $Q^{x,y}$.

\begin{lemma} \label{lem-4.1}
For all $\nu_{j}\in\mathcal{R}^{+} $, $j=1,\ldots,k$,
%
\begin{eqnarray}\label{nit8.2q0}
Q^{x,y}\Biggl( \prod_{j=1}^{k}L^{\nu_{j}}_{\infty}
\Biggr)&=&\sum_{\pi\in
\mathcal{P}_{k}} \int u(x, y_{1})u(y_{1},y_{2})
\nonumber
\\[-8pt]
\\[-8pt]
\nonumber
&&\hspace*{32pt}{} \cdots u(y_{k-1},y_{k })u(y_{k},y) \prod
_{j=1}^{k} \,d\nu_{\pi(j)}(y_{j}),
\nonumber
\end{eqnarray}
where the $\mathcal{P}_{k}$ denotes the set of permutations of $[1,k]$.
\end{lemma}

\begin{pf}
By (\ref{ls.3n}), we have
%
\begin{equation}
Q^{x,y}\Biggl( \prod_{j=1}^{k}L^{\nu_{j}}_{\infty}
\Biggr)=\int_{0}^{\infty} Q_{t}^{x,y}
\Biggl(\prod_{j=1}^{k}L^{\nu_{j}}_{\infty}
\circ k_{t} \Biggr) \,dt.\label{ls.3pq}
\end{equation}
Following the argument in (\ref{ls.1z}) and then using Lemma~\ref
{lem-finq}, we see that
\begin{eqnarray*}
\hspace*{-6pt}&&Q_{t}^{x,y}\Biggl( \Biggl(\prod
_{j=1}^{k}L^{\nu_{j}}_{\infty}\Biggr)
\circ k_{t}\Biggr)\\
\hspace*{-6pt}&&\qquad= \sum_{\pi
\in\mathcal{P}_{k}}
Q_{t}^{x,y}\biggl( \int_{0\leq r_{1}\leq\cdots\leq
r_{k}\leq t}
\,dL^{\nu_{\pi(1)}}_{r_{1}}\cdots dL^{\nu_{\pi
(k)}}_{r_{k}}\biggr)
\\
\hspace*{-6pt}&&\qquad= \sum_{\pi\in\mathcal{P}_{k}} \int_{0\leq r_{1}\leq\cdots\leq
r_{k}\leq t}\int
p_{r_{1}}(x,y_{1}) p_{r_{2}-r_{1}}(y_{1},y_{2})
\\
\hspace*{-6pt}&&\qquad\qquad\hspace*{82pt} \cdots p_{r_{k}-r_{k-1}}(y_{k-1},y_{k}) p_{t-r_{k}}(y_{k},y)
\prod_{j=1}^{k} \,d\nu_{\pi(j)}(y_{j})
\,dr_{j}.
\end{eqnarray*}
Therefore,
\begin{eqnarray*}
&&Q^{x,y}\Biggl( \prod_{j=1}^{k}L^{\nu_{j}}_{\infty}
\Biggr) \\
&&\qquad= \sum_{\pi\in\mathcal{P}_{k}} \int_{0\leq r_{1}\leq\cdots\leq
r_{k}\leq t<\infty}
\int p_{r_{1}}(x,y_{1}) p_{r_{2}-r_{1}}(y_{1},y_{2})
\\
&&\hspace*{107pt}\qquad\quad{}\cdots p_{r_{k}-r_{k-1}}(y_{k-1},y_{k}) p_{t-r_{k}}(y_{k},y)
\\
&&\hspace*{107pt}\qquad\quad{}\times \prod_{j=1}^{k} \,d\nu_{\pi(j)}(y_{j})
\,dr_{j} \,dt,
\end{eqnarray*}
which gives (\ref{nit8.2q0}).
\end{pf}

Let $\|\cdot\|$ be a proper norm. It follows from (\ref{nit8.2q0}) and
(\ref{8.2t}) that for any $\rho$ and $\nu_{1},\ldots,\nu_{k}\in
\mathcal{R}_{\|\cdot\|}$
%
\begin{equation}
\int Q^{x,x}\Biggl( \prod_{j=1}^{k}L^{\nu_{j}}_{\infty}
\Biggr) \,d\rho(x)=\mu\Biggl( L_{\infty
}^{\rho} \prod
_{j=1}^{k}L^{\nu_{j}}_{\infty}\Biggr).
\label{duh2j}
\end{equation}
For any $ \phi, \rho\in\mathcal{R}_{\|\cdot\|}^{+}$, set
%
\begin{equation}
Q^{\rho}_{\phi}(A)=\int Q^{x,x}\bigl(L^{\phi}_{\infty}1_{\{A\}}
\bigr) \,d\rho (x).\label{nit.1j}
\end{equation}
Note that by (\ref{duh2j}) we have $Q^{\rho}_{\phi}(\Om)=\mu
(L^{\rho
}_{\infty}L^{\phi}_{\infty})$, so that $Q^{\rho}_{\phi}$ is a
finite measure.
Using (\ref{duh2j}) again, we see that
%
\begin{equation}
Q^{\rho}_{\phi}\Biggl( \prod_{j=1}^{k}L^{\nu_{j}}_{\infty}
\Biggr)=\mu\Biggl( L_{\infty}^{\rho} L^{\phi}_{\infty}
\prod_{j=1}^{k}L^{\nu_{j}}_{\infty}
\Biggr) \label{duh3j}
\end{equation}
for all
$\nu_{j} \in\mathcal{R}_{\|\cdot\|} $.

By (\ref{8.2t}) and (\ref{duh3j}) and the fact that $\|\cdot\|$ is a
proper norm, we see that
%
\begin{equation}
\bigl|Q^{\rho}_{\phi}\bigl( \bigl( L_{\infty}^{\nu}
\bigr)^{n} \bigr)\bigr|=\bigl|\mu\bigl( L_{\infty
}^{\rho}
L^{\phi
}_{\infty}\bigl( L_{\infty}^{\nu}
\bigr)^{n} \bigr)\bigr| \leq n!C^{n}\|\phi\|\|\rho\|\|\nu
\|^{n}.\label{4.9}
\end{equation}
Therefore, $L^{\nu}_{\infty}$
is exponentially integrable with respect to the finite measures
$Q^{\rho
}_{\phi}(\cdot)$ and $\mu( L_{\infty}^{\rho} L_{\infty}^{\phi
}(\cdot) )$, so
that the finite dimensional distributions of $\{L^{\nu}_{\infty},
\nu\in
\mathcal{R}_{\|\cdot\|} \}$ under these measures are determined by
their moments. Consequently, by (\ref{duh3j}), for all bounded
measurable functions $F$ on $R^{k}$,
%
\begin{equation}
Q^{\rho}_{\phi}\bigl( F\bigl(L^{\nu_{1}}_{\infty},
\ldots, L^{\nu
_{k}}_{\infty}\bigr)\bigr)=\mu\bigl( L_{\infty}^{\rho}
L_{\infty}^{\phi} F\bigl(L^{\nu_{1}}_{\infty
},\ldots,
L^{\nu
_{k}}_{\infty}\bigr)\bigr). \label{palmers}
\end{equation}

We now obtain a Dynkin type isomorphism theorem that relates
permanental fields with kernel $u$ to continuous additive functionals
of a Markov process with
potential density $u$.
We can do this very efficiently by employing
a special case of the Palm formula for Poisson processes $\mathcal{L}$
with intensity measure $\xi$ on a measurable space $\mathcal{S}$, see
\cite{Bertoin}, Lemma~2.3, which states that for any positive function
$f$ on $\mathcal{S}$ and any measurable functional $G$ of $\mathcal{L}$
%
\begin{equation}
E_{\mathcal{L}} \biggl(\biggl(\sum_{\omega\in\mathcal{L}}f(\omega)
\biggr)G(\mathcal{L})\biggr)= \int E_{\mathcal{L}} \bigl(G\bigl(
\omega'\cup\mathcal{L}\bigr)\bigr)f\bigl(\omega'\bigr) \,d
\xi \bigl(\omega'\bigr). \label{palm.0}
\end{equation}

For $ \phi, \rho\in\mathcal{R}_{\|\cdot\|}^{+}$, we define
%
\begin{equation}
\th^{\rho,\phi} =\sum_{\omega\in\mathcal{L}_{\al}}L_{\infty}^{\rho}(
\omega) L_{\infty}^{\phi}(\omega ).\label{j8.5k}
\end{equation}
Obviously, $\th^{\rho,\phi}\ge0$.

\begin{theorem}[(Isomorphism Theorem I)] \label{theo-ljrm} Let $X$ be a transient
Borel right process with potential density $u$ as described in Section~\ref{ss-caf}. Let $\|\cd\|$ be a proper norm for $u$ and let $\{
\widetilde
\psi
(\nu),\nu\in\RR_{\|\cd\|}\}$ be as described in (\ref{ls.10a}). (By
Theorem~\ref{theo-9.1}, $\{\widetilde \psi(\nu),\nu\in\mathcal
{R}_{\|
\cdot\|
} \}$ is an $\al$-permanental field with kernel $u$.) Let
$\{L_{\infty}^{\nu},\nu\in\mathcal{R}_{\|\cdot\|} \}$ be as described
in the paragraph containing (\ref{potz.2}). Then
for any $ \phi, \rho\in\mathcal{R}_{\|\cdot\|}$ and all measures
$\nu
_{j}\in\mathcal{R}_{\|\cdot\|}$, $j=1,2,\ldots,$ and all bounded
measurable functions
$F$ on~$R^\infty$,
%
\begin{equation}
E_{\mathcal{L}_{\al}} Q_{\phi}^{\rho}\bigl( F \bigl( \widetilde
\psi(\nu _{i})+L_{\infty}^{\nu
_{i}} \bigr) \bigr)=
{1 \over\al}E_{\mathcal{L}_{\al}} \bigl(\th^{\rho,\phi} F\bigl(
\widetilde \psi (\nu_{i})\bigr)\bigr),\label{6.1}
\end{equation}
and $\th^{\rho,\phi} $, given in (\ref{j8.5k}), has all its moments finite.
[Here, we use the notation $F( f(x_{ i})):=F( f(x_{ 1}),f(x_{
2}),\ldots
)$.]
\end{theorem}

Since (\ref{6.1}) depends only on the distribution of the $\al
$-permanental field with kernel $u$, this theorem implies Theorem~\ref
{theo-ljrmi}.

\begin{pf}
We apply the Palm formula with intensity measure $\al\mu$,
%
\begin{equation}
f(\omega)= L_{\infty}^{\rho}(\omega) L_{\infty}^{\phi
}(
\omega)\label{palm.2}
\end{equation}
and
%
\begin{equation}
G(\mathcal{L}_{\al})=F\bigl( \widehat L_{\delta}^{\nu_{j}}
\bigr).\label{palm.2a}
\end{equation}
To begin let $F$ be a bounded continuous function on $R^{n}$. Note that
%
\begin{equation}
\sum_{\omega\in\mathcal{L}_{\al}}f(\omega)=\th^{\rho,\phi
}.\label{palm.2aa}
\end{equation}
Note also that since $\omega'$ and $\LL_{\al}$ are disjoint a.s.
%
\begin{eqnarray}
\widehat L_{\delta}^{\nu_{j}}\bigl(\omega'\cup
\mathcal{L}_{\al}\bigr) &=&\biggl(\sum_{\omega\in\omega'\cup\mathcal{L}_{\al}}
1_{\{\ze
(\omega)>\delta\}} L_{\infty
}^{\nu_{j}}(\omega)\biggr)-\al\mu\bigl(
1_{\{\ze>\delta\}}L_{\infty
}^{\nu_{j}}\bigr)
\nonumber
\\[-8pt]
\\[-8pt]
\nonumber
&=&1_{\{\ze(\omega')>\delta\}} L_{\infty}^{\nu_{j}}\bigl(\omega '
\bigr)+\widehat L_{\delta}^{\nu
_{j}}(\mathcal{L}_{\al}),\label{palm.2b}
\end{eqnarray}
so that
%
\begin{equation}
G\bigl(\omega'\cup\mathcal{L}_{\al}\bigr)=F\bigl(\widehat
L_{\delta}^{\nu
_{j}}(\mathcal{L}_{\al
})+1_{\{\ze(\omega')>\delta\}}
L_{\infty}^{\nu_{j}}\bigl(\omega '\bigr)
\bigr).\label{palm.3}
\end{equation}
It follows from (\ref{palm.0}) that
%
\begin{eqnarray}\label{palm.4}\qquad
&& E_{\mathcal{L}_{\al}} \bigl(\th^{\rho,\phi} F\bigl(\widehat
L_{\delta
}^{\nu
_{j}}\bigr)\bigr)
\nonumber
\\[-8pt]
\\[-8pt]
\nonumber
&&\qquad =\al\int E_{\mathcal{L}_{\al}} \bigl(L_{\infty}^{\rho}\bigl(
\omega'\bigr) L_{\infty}^{\phi}\bigl(\omega
'\bigr) F\bigl(\widehat L_{\delta}^{\nu_{j}}(
\mathcal{L}_{\al})+1_{\{\ze
(\omega
')>\delta\}} L_{\infty}^{\nu_{j}}
\bigl(\omega'\bigr)\bigr) \bigr) \,d\mu\bigl(\omega'
\bigr).
\nonumber
\end{eqnarray}
We interchange the integrals on the right-hand side of (\ref{palm.4})
and use (\ref{palmers}) and then take the limit as $\delta\rar0$,
to get
(\ref{6.1})
for bounded continuous functions $F$ on $R^{n}$. The extension to
general bounded
measurable functions
$F$ on $R^\infty_+$ is routine.

To see that $\th^{\rho,\phi} $ has all moments finite, we use
the master formula for Poisson processes in the form
%
\begin{equation}
E_{\mathcal{L}_{\al}}\bigl(e^{ z\th^{\rho,\phi}}\bigr)=\exp\biggl(\al\biggl(\int
_{\Om}\bigl( e^{z L^{\rho}_{\infty} L^{\th}_{\infty} } -1\bigr) \,d\mu(\omega) \biggr)
\biggr) \label{ls.11t}
\end{equation}
with $z<0$.
Differentiating each side of (\ref{ls.11t}) $n$ times with respect to
$z$ and then taking $ z\uparrow0$ we see that
%
\begin{equation}
E_{\mathcal{L}_{\al}}\bigl( \bigl(\th^{\rho,\phi}\bigr)^{n} \bigr)=
\sum_{\bigcup_{i }
B_{i}=[1,n]} \prod_{i }
\al\mu\bigl(\bigl(L^{\rho}_{\infty} L^{\th
}_{\infty}
\bigr)^{|
B_{i}|}\bigr),\label{mom.11t}
\end{equation}
where the sum is over all partitions $B_{1},\ldots, B_{n}$ of $[1,n]$.
This is finite for $ \phi, \rho\in\mathcal{R}_{\|\cdot\|}^{+}$.
\end{pf}

Isomorphism Theorem I shows that the continuity of the permanental
field implies the continuity (in the measures) of the continuous
additive functionals.

\begin{corollary}\label{cor-4.1} In the notation and under the
hypotheses of Theorem~\ref{theo-ljrm}, let $D\subseteq\RR_{\|\cd\|
} $
and suppose there exists a metric $d$ on $D$ such that
%
\begin{equation}
\lim_{\delta\rar0} E_{\mathcal{L}_{\al}} \Bigl(\mathop{\sup_{d(\nu
,\nu
')\le\delta}}_{\nu,\nu'\in D}
\bigl|\psi(\nu)-\psi\bigl(\nu'\bigr) \bigr|^2\Bigr)=0,\label{6p}
\end{equation}
then
%
\begin{equation}
\lim_{\delta\rar0}Q^{\rho}_{\phi} \Bigl(\mathop{\sup
_{d(\nu,\nu
')\le\delta}}_{
{\nu,\nu'\in D} } \bigl|L^{\nu}_{\infty} -L^{\nu'}_{\infty}
\bigr|\Bigr)=0.\label{nit.6po}
\end{equation}
\end{corollary}

\begin{pf}
It
follows from (\ref{6.1}) that
%
\begin{eqnarray}\label{12.1}
&& Q^{\rho}_{\phi}\Bigl( \mathop{\sup_{d(\nu,\nu')\le\delta}}_ {\nu
,\nu'\in D}
\bigl|L^{\nu}_{\infty}-L^{\nu'}_{\infty}\bigr|\Bigr)
\nonumber
\\
&&\qquad \le E_{\mathcal{L}_{\al}} \Bigl(\mathop{\sup_{d(\nu,\nu')\le
\delta}}_{\nu
,\nu'\in D } \bigl| \psi(\nu)-\psi
\bigl(\nu'\bigr) \bigr| \Bigr) Q^{\rho
}_{\phi
}(1)
\\
\nonumber
&&\qquad\quad{} + {1 \over\al} \Bigl( E_{\mathcal{L}_{\al}} \Bigl(\mathop{\sup
_{d(\nu,\nu')\le\delta}}_ {\nu,\nu'\in D} \bigl|\psi(\nu)-\psi\bigl(\nu'\bigr)
\bigr|^2 \Bigr) E_{\mathcal{L}_{\al
}} \bigl(\th ^{\rho,\phi}
\bigr)^{2} \Bigr)^{1/2}.
\nonumber
\end{eqnarray}
Using this, it is easy to see that (\ref{6p}) implies (\ref{nit.6po}).
\end{pf}

For applications of the isomorphism theorem in Section~\ref{sec-caf},
we sometimes need to consider measures $\rho$ and $\phi$ in (\ref{6.1})
that are not necessarily in $ \mathcal{R}_{\|\cdot\|} $.
To deal with this, we introduce two additional norms on $\MM(S)$:
%
\begin{equation}
\|\nu\|_{u^2,\infty}:= |\nu|(S)\vee\sup_{x}\int\bigl(
u^{2}(x,y) + u^{2}(y,x) \bigr) \,d|\nu| (y),\label{p.2}
\end{equation}
where $|\nu|$ is the total variation of the measure $\nu$, and
%
\begin{equation}
\|\nu\|_{0}:=|\nu|(S)\vee\sup_{x}\int u(x,y) \,d|
\nu| (y).\label{0norm}
\end{equation}

\begin{lemma}\label{lem-prop2}
Let $A\cup B$ be a partition of $[1,n]$, $n\geq2$, with $B\neq
\varnothing$. Then
%
\begin{eqnarray}\label{bp.5}
&& \Biggl|\int u (y_{1},y_{2})\cdots u (y_{n-1},y_{n})u
(y_{n},y_{1}) \prod_{j=1}^{n}
\,d\nu_{j}(y_{j}) \Biggr|
\nonumber
\\[-8pt]
\\[-8pt]
\nonumber
&& \qquad\leq\prod_{i\in A}\|\nu_{i}
\|_{0}\prod_{j\in B} \|\nu_{j}\|
_{u^2,\infty}.
\nonumber
\end{eqnarray}
Let $ \phi\in\mathcal{R}_{\|\cdot\|_{u^2,\infty}}^{+}$ and $\rho
\in
\mathcal{R}_{\|\cdot\|_{0}}^{+}$. In addition, let $\nu_{j}\in
\mathcal
{R}_{\|\cdot\|}, j=1,\ldots,k $, for some proper norm $\|\cdot\|$. Then
there exists a constant $C=C( \phi,\rho, \|\cdot\|)<\infty$, such that
%
\begin{equation}
\Biggl|\mu\Biggl( L_{\infty}^{\rho} L^{\phi}_{\infty}\prod
_{j=1}^{k}L^{\nu
_{j}}_{\infty}
\Biggr) \Biggr|\leq k!C^{k} \|\rho\|_{0} \|\phi\| _{u^2,\infty}
\prod_{j=1}^{k}\|\nu_{j}\|.
\label{bp.5m}
\end{equation}
\end{lemma}

\begin{pf} Without loss of generality, we assume that $1\in B$. Then
using the Cauchy--Schwarz inequality in $y_{1}$, we have
%
\begin{eqnarray}\label{jd.10a}
&&\Biggl|\int u (y_{1},y_{2})\cdots u
(y_{k-1},y_{k})u (y_{k},y_{1}) \prod
_{j=1}^{k} \,d\nu_{j}(y_{j})
\Biggr|
\\
&&\qquad \leq\int\biggl(\int u^{2} (y_{1},y_{2}) \,d|
\nu_{1}| (y_{1})\biggr)^{1/2} \biggl(\int
u^{2} (y_{k},y_{1}) \,d|\nu_{1}|
(y_{1})\biggr)^{1/2}
\nonumber
\\
&&\qquad\qquad \cdots u (y_{2},y_{3})\cdots u (y_{k-1},y_{k})
\prod_{j=2}^{k} \,d|\nu_{j}|(y_{j})
\nonumber
\\
&&\qquad\leq\|\nu_{1}\|_{u^2,\infty}\int u(y_{2},y_{3})
\cdots u (y_{k-1},y_{k}) \prod_{j=2}^{k}
\,d|\nu_{j}|(y_{j}).
\nonumber
\end{eqnarray}
We bound successively the integrals with respect to $d|\nu_{j}|(y_{j})$
for $j=k,k-1,\ldots, 3$
to obtain
%
\begin{eqnarray}\label{jd.11}
&& \int u (y_{2},y_{3})\cdots u (y_{k-1},y_{k})
\prod_{j=2}^{k} \,d|\nu _{j}|(y_{j}).
\nonumber\\
&&\qquad \leq\biggl(\sup_{x}\int u(x,y) \,d|\nu_{k}| (y)
\biggr)\int u (y_{2},y_{3})\cdots u (y_{k-2},y_{k-1})
\prod_{j=2}^{k-1} \,d|\nu _{j}|(y_{j})
\\
&& \qquad\leq\prod_{j=3}^{k} \biggl(\sup
_{x}\int u(x,y) \,d|\nu_{j}| (y)\biggr)\int1 \,d|
\nu_{2}|(y_{2}) \le\prod_{j=2}^{k}
\|\nu_{j}\|_{0}.
\nonumber
\end{eqnarray}
Using (\ref{jd.10a}) and (\ref{jd.11}) we see that
%
\begin{eqnarray}\label{jd.11q}
&&\Biggl|\int u (y_{1},y_{2})\cdots u (y_{n-1},y_{n})u
(y_{n},y_{1}) \prod_{j=1}^{n}
\,d\nu_{j}(y_{j}) \Biggr|
\nonumber
\\[-8pt]
\\[-8pt]
\nonumber
&&\qquad\leq\|\nu_{1}
\|_{u^2,\infty
}\prod_{j=2}^{n}\|
\nu_{j}\|_{0}.
\end{eqnarray}

We now note that by the Cauchy--Schwarz inequality for the finite
measure $\nu$, $\|\nu\|_{0}\leq C\|\nu\|_{u^2,\infty}$. Using this and
(\ref{jd.11q}) and recognizing that the choice of indices in (\ref
{jd.11q}) is arbitrary, we get (\ref{bp.5}).

To obtain (\ref{bp.5m}), we use the
Cauchy--Schwarz inequality to write
%
\begin{equation}\qquad
\Biggl|\mu \Biggl( L_{\infty}^{\rho} L^{\phi}_{\infty}
\prod_{j=1}^{k}L^{\nu
_{j}}_{\infty}
\Biggr) \Biggr|\leq \bigl\{ \mu\bigl(\bigl( L_{\infty}^{\rho}
L^{\phi}_{\infty
}\bigr)^{2}\bigr) \bigr\}
^{1/2} \Biggl\{ \mu \Biggl( \Biggl( \prod_{j=1}^{k}L^{\nu
_{j}}_{\infty
}
\Biggr)^{2} \Biggr) \Biggr\}^{1/2}. \label{duh3a}
\end{equation}
We use (\ref{bp.5}) with $|A|=|B|=2$ to bound the first term and (\ref
{8.2t}) and (\ref{p.0}), and the fact that $((2k-1)!)^{1/2}\le C^{k}
k!$ to bound the second term and get (\ref{bp.5m}).
\end{pf}

Using Lemma~\ref{lem-prop2}, we can modify the hypotheses of Theorem~\ref{theo-ljrm} to obtain a second isomorphism theorem.

\begin{theorem}[(Isomorphism Theorem II)] \label{theo-ljr} All the
results of Theorem~\ref{theo-ljrm} hold for $ \phi\in\mathcal{R}_{\|\cdot\|
_{u^2,\infty
}}^{+} $ and $\rho\in\mathcal{R}_{\|\cdot\|_{0}}^{+}$.
\end{theorem}

\begin{pf} Given the proof of Theorem~\ref{theo-ljrm}, to prove this
theorem it suffices to show that (\ref{palmers}) holds when $ \phi\in
\mathcal{R}_{\|\cdot\|_{u^2,\infty}}^{+} $ and $\rho\in\mathcal
{R}^{+}_{\|\cdot\|_{0}}$. To do this, we first show that the argument
from (\ref{duh2j})--(\ref{palmers}) holds under this change of hypothesis.

Set
%
\begin{equation}
Q^{\rho}_{\phi}(A)=\int Q^{x,x}\bigl(L^{\phi}_{\infty}1_{\{A\}}
\bigr) \,d\rho (x).\label{nit.1}
\end{equation}
By Remark~\ref{rem-3.1}, (\ref{8.2t}) holds for measures in $\RR^{+}$.
In particular, by (\ref{nit8.2q0}), for $\rho, \nu_{j}\in\mathcal{R}^{+}$,
%
\begin{equation}
\int Q^{x,x}\Biggl( \prod_{j=1}^{k}L^{\nu_{j}}_{\infty}
\Biggr) \,d\rho(x)=\mu\Biggl( L_{\infty
}^{\rho} \prod
_{j=1}^{k}L^{\nu_{j}}_{\infty}\Biggr).
\label{duh2}
\end{equation}
Therefore, $Q^{\rho}_{\phi}(\Om)=\mu(L^{\rho}_{\infty}L^{\phi
}_{\infty})$ so
that by Lemma~\ref{lem-prop2}, (\ref{bp.5}), we see that $Q^{\rho
}_{\phi
}$ is a finite measure. Using (\ref{duh2}), we see that
%
\begin{equation}
Q^{\rho}_{\phi}\Biggl( \prod_{j=1}^{k}L^{\nu_{j}}_{\infty}
\Biggr)=\mu\Biggl( L_{\infty}^{\rho} L^{\phi}_{\infty}
\prod_{j=1}^{k}L^{\nu_{j}}_{\infty}
\Biggr) \label{duh3}
\end{equation}
for all $\nu_{j} \in\mathcal{R}^{+} $.

We now use Lemma~\ref{lem-prop2}, (\ref{bp.5m}), to see that (\ref
{duh3}) holds for $ \phi\in\mathcal{R}_{\|\cdot\|_{u^2,\infty
}}^{+} $,
$\rho\in\RR_{\|\cdot\|_{0}}^{+}$ and $\{\nu_{j}\}\in\mathcal
{R}_{\|
\cdot\|}$. Therefore, using Lemma~\ref{lem-prop2} we see that for
any $\nu\in\mathcal{R}_{\|\cdot\|} $
%
\begin{equation}
\bigl|Q^{\rho}_{\phi}\bigl( \bigl( L_{\infty}^{\nu}
\bigr)^{n} \bigr)\bigr|=\mu\bigl( L_{\infty
}^{\rho}
L^{\phi
}_{\infty}\bigl( L_{\infty}^{\nu}
\bigr)^{n} \bigr) \leq n!C^{n} \|\rho\|_{0} \|
\phi\|_{u^2,\infty}\|\nu\|^{n}, \label{duh3aa}
\end{equation}
which shows that all $\{L^{\nu_{j}}_{\infty}\}$ are exponentially
integrable with respect to the finite measures $Q^{\rho}_{\phi}$ and
$\mu( L_{\infty}^{\rho} L_{\infty}^{\phi}\cd)$. Since (\ref
{duh3}) holds for
$ \phi\in\mathcal{R}_{\|\cdot\|_{u^2,\infty}}^{+} $, $\rho\in
\RR_{\|\cdot
\|_{0}}^{+}$ and $\{\nu_{j}\}\in\mathcal{R}_{\|\cdot\|}$; this
shows that
for all bounded measurable functions $F$ on $R^{k}$
%
\begin{equation}
Q^{\rho}_{\phi}\bigl( F\bigl(L^{\nu_{1}}_{\infty},
\ldots, L^{\nu
_{k}}_{\infty}\bigr)\bigr)=\mu\bigl( L_{\infty}^{\rho}
L_{\infty}^{\phi} F\bigl(L^{\nu_{1}}_{\infty
},\ldots,
L^{\nu
_{k}}_{\infty}\bigr)\bigr) \label{palmersa}
\end{equation}
holds when $ \phi\in\mathcal{R}_{\|\cdot\|_{u^2,\infty}}^{+} $,
$\rho\in
\RR_{\|\cdot\|_{0}}^{+}$ and $\{\nu_{j}\}\in\mathcal{R}_{\|\cdot\|}$.
With this modification, the proof of Theorem~\ref{theo-ljrm} proves
this theorem.
\end{pf}

\begin{remark}\label{rem-4.1} It is easy to see that Corollary~\ref
{cor-4.1} also holds under the hypotheses of Theorem~\ref{theo-ljr}.
\end{remark}

\section{Continuity of \texorpdfstring{$L^{\nu}_{\infty}$}{$L {nu} {infty}$} and \texorpdfstring{$\psi(\nu)$}{{$psi(nu)$}}}\label
{sec-contperm}

In this section, we give sufficient conditions for the continuity of
the additive functionals $\{L^{\nu}_{\infty}, \nu\in\VV\}$ and
permanental fields $\{\psi(\nu), \nu\in\VV\}$ that extend well-known
results for second-order Gaussian chaoses.

Let $(T,\tau)$ be a~metric or pseudo-metric space. Let $B_{\tau
}(t,u)$ denote the closed ball in $(T,\tau)$ with radius $u$ and center
$t$. For any probability measure $\si$ on $(T,\tau)$, we define
%
\begin{equation}
J_{T, \tau,\si}( a) =\sup_{t\in T}\int_0^a
\log\frac{1}{\si
(B_{\tau
}(t,u))} \,du.\label{tau}
\end{equation}

Let $\VV$ be a linear space of measures on $S$ and $u$ a kernel on
$S\times S$. Suppose that $\|\cdot\|$ is a proper norm for $\VV$ with
respect to $u$ and $\VV\subseteq\mathcal{R}_{\|\cdot\|} $. Then $\|
\nu
-\nu'\|$ is a metric on $\VV$. In this situation, we write $J_{T,
\tau
,\si}$ in (\ref{tau}) as $J_{\VV, \|\cdot\|,\si}$.

\begin{theorem}\label{cor-jack}
Let $ \phi\in\mathcal{R}_{\|\cdot\|_{u^2,\infty}}^{+}$ and $\rho
\in
\mathcal{R}_{\|\cdot\|_{0}}^{+}$, and let $\|\cdot\|$ be a proper norm.
Assume that there exists
a probability measure $\si$ on $\VV$ such that $J_{\VV, \|\cdot\|
,\si
}(D)<\infty$, where $D$ is the diameter of $\VV$ with respect to $ \|
\cdot
\|$ and
%
\begin{equation}
\lim_{\delta\to0} J_{\VV, \|\cdot\|,\si}(\delta)=0. \label{2.10vi}
\end{equation}
Then for any countable set $D\subseteq\VV$, with compact closure
%
\begin{equation}
\lim_{\delta\rar0}Q^{\rho}_{\phi} \Bigl(\mathop{\sup
_{d(\nu,\nu
')\le\delta}}_
{\nu,\nu'\in D}  \bigl|L^{\nu}_{\infty} -L^{\nu'}_{\infty}
\bigr|\Bigr)=0.\label{nit.6pjack}
\end{equation}

A similar result holds for $\{\psi(\nu),\nu\in\VV\}$ with respect
to $
E_{\mathcal{L}_{\al}}$.
\end{theorem}

\begin{pf*}{Proof of Theorems \ref{theo-contprop} and \ref
{cor-jack}} These theorems are immediate consequence of the following
lemma and the well-known sufficient condition for continuity of
stochastic processes with a metric in an exponential Orlicz space; see,
for example, \cite{MR09}, Section~3, or \cite{MRpermfields}, Theorem~2.1.
\end{pf*}

Let $
\Xi(x)=
\exp(x ) -1
$ and $L^{\Xi}(\Om,\FF,P)$ denote the set of random
variables
$\xi\dvtx \Om\to R^{1}$ such that $E(\Xi( |\xi|/c ))<\infty$ for some $c>0$.
$L^{\Xi}(\Om,\FF,P)$ is a Banach space with norm given by
%
\begin{equation}
\|\xi\|_{\Xi}=\inf \bigl\{c>0\dvtx E\bigl(\Xi\bigl( |\xi|/c\bigr)\bigr)\le1 \bigr
\}.\label{1.4}
\end{equation}

\begin{lemma}\label{cor-2.1}Let $ \phi\in\mathcal{R}_{\|\cdot\|
_{u^2,\infty
}}^{+}$ and $\rho\in\mathcal{R}_{\|\cdot\|_{0}}^{+}$, and let $\|
\cdot
\|$ be a proper norm. Then there exists a constant $C=C( \phi,\rho, \|
\cdot\|)<\infty$, such that
%
\begin{equation}
\bigl\| L^{\nu}_{\infty}\bigr\|_{\Xi}\le C \| \nu\|\qquad \forall\nu\in
\VV, \label{1.9dj}
\end{equation}
where $ \| \cd\|_{\Xi}$ is the norm of the exponential Orlicz space
generated by $e^{|x|}-1$ with respect to $Q^{\rho}_{\phi}$.

Similarly, let $\{\psi(\nu),\nu\in\VV\}$ be an $\al$-permanental field
with kernel $u$ and $\|\cdot\|$ a proper norm with respect to $u$, then
for some $C_{\al}<\infty$, depending only on $\al$,
%
\begin{equation}
\bigl\| \psi(\nu)\bigr\|_{\Xi}\le C_{\al} \| \nu\|\qquad \forall\nu\in\VV,
\label{1.9d}
\end{equation}
where $ \| \cd\|_{\Xi}$ is the norm of the exponential Orlicz space
generated by $e^{|x|}-1$ with respect to $ E_{\mathcal{L}_{\al}}$.
\end{lemma}

\begin{pf} Since $Q^{\rho}_{\phi}$ is a finite measure, it follows from
the Cauchy--Schwarz inequality and (\ref{duh3aa}) that
%
\begin{equation}
Q^{\rho}_{\phi}\bigl(\bigl| L^{\nu}_{\infty}\bigr|^{n}
\bigr) \leq C \bigl(Q^{\rho}_{\phi
}\bigl( \bigl(L^{\nu
}_{\infty}
\bigr)^{2n}\bigr)\bigr)^{1/2}\leq n!C^{n}\|\nu
\|^{n}.\label{qeasy1}
\end{equation}
The inequality in (\ref{1.9dj}) follows from this.

The inequality in (\ref{1.9d}) can be derived similarly, using the
Cauchy--Schwarz inequality, Definition~\ref{def1}, the definition of
proper norms (\ref{p.0}), and the fact that there are $n!$ permutations
of $[1,n]$.
\end{pf}

Other results on the continuity of permanental fields are given in
\cite{MRpermfields}.

\begin{remark}\label{rem-4.1a} Using the isomorphism theorem,
(\ref{nit.6pjack}) can be derived from the similar result for $\{\psi
(\nu),\nu\in\VV\}$; see Remark~\ref{rem-4.1}. In all our earlier work,
continuity conditions for local times and other continuous additive
functionals of Markov processes are obtained in this way, that is, by
means of an isomorphism theorem. It is noteworthy that in this paper
(\ref{nit.6pjack}) is obtained directly using properties of the loop measure.
\end{remark}

\section{Joint continuity of continuous additive functionals} \label{sec-caf}

In this section, we obtain sufficient conditions for continuity
of the stochastic process
%
\begin{equation}
L= \bigl\{ L^{\nu}_{t}, (t,\nu)\in R^{+}\times\VV
\bigr\} \label{sum}
\end{equation}
for some family of measures $\VV\subseteq\mathcal{R}^{+}$, endowed
with a topology induced by an appropriate proper norm.

By definition, $L^{\nu}_{t}$ is continuous in $t$. However, proving
the joint continuity of~(\ref{sum}), $P^{x}$ almost surely, is
difficult. We break the proof into a series of lemmas and theorems.
We assume that $ \phi\in\mathcal{R}_{\|\cdot\|_{u^2,\infty}}^{+}$ and
$\rho\in\mathcal{R}_{\|\cdot\|_{0}}^{+}$, which as noted above,
implies that
$Q^{\rho}_{\phi}$ [defined in (\ref{nit.1})], is a finite measure.

Let $h(x,y)$ be a bounded measurable function on S which is excessive
in $x$, and such that
%
\begin{equation}
0<h(x,y)\leq u(x,y),\qquad x,y\in S.\label{low.1}
\end{equation}
For example, we can take $h(x,y)=1\wedge u(x,y)$, or more generally
$h(x,y)=f(x)\wedge u(x,y)$ for any bounded strictly positive excessive
function $f$. In the proof of Theorem~\ref{theo-1.1}, we take $h(x,y)=
u_{1}(x,y)=\int_{1}^{\infty}p_{t}(x,y) \,dt$.

Set $h_{y}(z)=h(z,y)$. We let $Q^{x,h_{y}}$ denote the (finite) measure
defined by
%
\begin{equation}
Q^{x,h_{y}}( 1_{\{\ze>s \}} F_s)= P^x\bigl(
F_s h(X_s,y )\bigr)\qquad \mbox{for all $F_s
\in b\mathcal{F}^{0}_s$}, \label{nit4.4a}
\end{equation}
where $\mathcal{F}^{0}_s$ is the $\sigma$-algebra generated by $\{
X_{r}, 0\leq r\leq s\}$. In this notation, we can write the $\si
$-finite measure $Q^{x,y}$ in (\ref{nit4.4}) as $Q^{x,u_{y}}$.

Set
%
\begin{equation}
Q^{x,h_{x}}_{\phi}(A)=Q^{x,h_{x}}\bigl(L^{\phi}_{\infty}1_{\{A\}}
\bigr)\label{nit.1l}
\end{equation}
and
%
\begin{equation}
Q^{\rho, h}_{\phi}(A)=\int Q^{x,h_{x}}\bigl(L^{\phi}_{\infty}1_{\{A\}}
\bigr) \,d\rho (x)=\int Q^{x,h_{x}}_{\phi}(A) \,d\rho(x).\label{nit.1+}
\end{equation}
Note that it follows from (\ref{low.1}) and (\ref{nit4.4a}) that
%
\begin{equation}
Q^{\rho, h}_{\phi}(A)\leq Q^{\rho}_{\phi}(A)\label{low.2}
\end{equation}
for all $A\in\mathcal{F}^{0}$.

\begin{lemma}\label{nit1}Let $X =
(\Om, X_t,P^x
)$ be a Borel right process in $S$ with strictly positive
potential densities
$u(x,y)$, and let
$\VV\subseteq\RR^{+}_{\|\cdot\|}$, where $\|\cdot\|$ is proper for
$u$. Let $\OO$ be a topology for $\VV$
under which $\VV$ is a separable locally compact metric space with
metric $d$. Assume that there exist
measures $\rho\in\mathcal{R}_{\|\cdot\|_{0}}^{+}$, and $ \phi\in
\mathcal{R}_{\|\cdot\|_{u^2,\infty}}^{+}$ for which:
\begin{longlist}[(ii)]
\item[{(i)}]
%
\begin{equation}\quad
\int u(y,z)h_{x}(z) \,d\nu(z)\quad\mbox{and}\quad \int u(y,z)\int
u(z,w)h_{x}(w) \,d\phi(w) \,d\nu(z)
\end{equation}
are continuous
in $\nu\in\VV$, uniformly in $y,x\in S$, and
\item[{(ii)}] for any countable set $D\subseteq\VV$, with
compact closure
%
\begin{equation}
\lim_{\delta\rar0}Q^{\rho}_{\phi} \Bigl(\mathop{\sup_{d(\nu,\nu
')\le\delta}}_
{\nu,\nu'\in D} \bigl|L^{\nu}_{\infty} -L^{\nu'}_{\infty}
\bigr|\Bigr)=0,\label{nit.6p}
\end{equation}
where $\{L_{t}^{\nu},\nu\in D\}$ are continuous additive functionals of
$X$ as defined in Section~\ref{ss-caf}.
\end{longlist}

Then
for any $\ep>0$, there exists a $\delta>0$, such that
%
\begin{equation}
Q^{\rho,h}_{\phi} \Bigl( \sup_{t\geq0}\mathop{\sup
_{d(\nu,\nu')\le
\delta}}_
{\nu,\nu'\in D}  L^{\nu}_{t}-L^{\nu'}_{t}
\ge2\ep\Bigr) \leq\ep.\label{p.1}
\end{equation}
\end{lemma}

\begin{pf} As in \cite{MR96}, we use martingale techniques to go from
(\ref{nit.6p}) to (\ref{p.1}). However, the present situation is
considerably more complicated.

By working locally it suffices to consider $ \VV$ compact.
For fixed $y$, let
%
\begin{equation}
P^{x/h_{y}}( \cdot)= {Q^{x,h_{y}}(\cdot) \over h_{y}(x)},\qquad x\in S.\label{low.3}
\end{equation}
$
(\Om, X_t,P^{x/h_{y}}
)$ is a Borel right process in $S$, called the $h_{y}$-transform of\break $
(\Om, X_t, P^x
)$, \cite{S}, Section~62.

To begin we fix $x\in S$. Set
%
\begin{equation}
Z={ L^{\phi}_{\infty}\over E^{x/h_{x}}(L^{\phi}_{\infty})}\quad \mbox{and}\quad Z_{s}=E^{x/h_{x}}\bigl(Z
\mid\mathcal F^{0}_s\bigr),\label{de.1a}
\end{equation}
and define the probability
measure
%
\begin{equation}
P_{\phi}^{x/h_{x}}(A ):= E^{x/h_{x}}(1_{\{A\}}Z)=
{ Q_{\phi
}^{x,h_{x}}(A ) \over Q^{x,h_{x}}(L^{\phi}_{\infty})}.\label{de.1}
\end{equation}

By \cite{book}, Lemma~3.9.1, we can assume that the continuous
additive functionals $L^{\nu
}_{t}$ are $\mathcal F^{0}_t$ measurable.
Consider the $ P_{\phi}^{x/h_{x}}$ martingale
%
\begin{equation}
A^{\nu}_s=E_{\phi}^{x/h_{x}}
\bigl(L^{\nu}_{\infty}\mid\mathcal F^{0}_s
\bigr)={E^{x/h_{x}}(L^{\nu}_{\infty}Z\mid\mathcal F^{0}_s) \over
Z_{s}}.\label{martj}
\end{equation}
The last equality is well known and easy to check. Using the additivity
property,
%
\begin{equation}
L^{\nu}_{\infty} =L^{\nu}_s+
L^{\nu}_{\infty}\circ\tau_s,\label{add}
\end{equation}
where $\tau_{s}$ denotes the shift operator on $\Om$, we see that
%
\begin{equation}
A^{\nu}_s =L^{\nu}_s+
{E^{x/h_{x}}(L^{\nu}_{\infty}\circ\tau_s
Z\mid
\mathcal F^{0}_s) \over Z_{s}}:=L^{\nu}_s+ H^{\nu}_s.\label{martjq}
\end{equation}
Let $D$ be a countable dense subset of $\VV$.
By (\ref{martjq}), for any finite
subset $F\subset D$,
%
\begin{eqnarray} \label{12.2a}
&& P_{\phi}^{x/h_{x}} \Bigl( \sup_{t\geq0}\mathop{\sup
_{d(\nu,\nu')\le
\delta
}}_ {\nu,\nu'\in F}  L^{\nu}_{t}-L^{\nu'}_{t}
\ge3\ep\Bigr)
\nonumber\\
&&\qquad \le P_{\phi}^{x/h_{x}} \Bigl(\sup_{t\geq0} \mathop{\sup
_{d(\nu
,\nu
')\le\delta}}_ {\nu,\nu'\in F} A^\nu_t-A^{\nu'}_t
\ge\ep\Bigr)
\nonumber
\\[-8pt]
\\[-8pt]
\nonumber
&&\qquad\quad {}+ P_{\phi}^{x/h_{x}} \Bigl(\sup_{t\geq0} \mathop{\sup
_{d(\nu,\nu
')\le\delta}}_ {\nu,\nu'\in F} H^\nu_t-H^{\nu'}_t
\ge2\ep\Bigr)
\\
&& \qquad:=I_{1,x} +I_{2,x}.
\nonumber
\end{eqnarray}

Using (\ref{add}), but this time for $L_{\infty}^{\phi}$, and using the
Markov property, we see that
%
\begin{eqnarray}\label{12.2as}
H^{\nu}_t&=&{E^{x/h_{x}}(L^{\nu}_{\infty}\circ\tau_t L^{\phi
}_{\infty}\mid
\mathcal F^{0}_t) \over E^{x/h_{x}}(L^{\phi}_{\infty}\mid\mathcal F^{0}_t)}
\nonumber\\
& = &{L^{\phi}_{t} E^{x/h_{x}}(L^{\nu}_{\infty}\circ\tau_t \mid
\mathcal
F^{0}_t)+E^{x/h_{x}}((L^{\nu}_{\infty}L^{\phi}_{\infty})\circ
\tau_t \mid
\mathcal F^{0}_t) \over E^{x/h_{x}}(L^{\phi}_{\infty}\mid\mathcal
F^{0}_t)}
\\
& =& {L^{\phi}_{t} E^{X_{t}/h_{x}}(L^{\nu}_{\infty
})+E^{X_{t}/h_{x}}(L^{\nu}_{\infty}L^{\phi}_{\infty} ) \over
E^{x/h_{x}}(L^{\phi}_{\infty}\mid\mathcal F^{0}_t)}.
\nonumber
\end{eqnarray}
Here and throughout, we are using the convention that $f(X_{t})=1_{\{
t>\ze\}}f(X_{t})$ for any function $f$ on $S$.
Proceeding the same way with the denominator, we obtain
%
\begin{equation}
H^{\nu}_t= {L^{\phi}_{t} E^{X_{t}/h_{x}}(L^{\nu}_{\infty
})+E^{X_{t}/h_{x}}(L^{\nu}_{\infty}L^{\phi}_{\infty} ) \over
L^{\phi
}_{t}+E^{X_{t}/h_{x}}(L^{\phi}_{\infty} )}.\label{12.2at}
\end{equation}
Using
\cite{MR96}, (2.25) and (2.22), where $u\bb(\cd)=h_{x}(\cd)$, we have
%
\begin{equation}
E^{y/h_{x}}\bigl(L^{\nu}_{\infty}\bigr)={\int u(y,z)h_{x}(z) \,d\nu(z) \over
h_{x}(y)}
\label{mema}
\end{equation}
and
%
\begin{eqnarray} \label{memb}
&&E^{y/h_{x}}\bigl(L^{\nu}_{\infty}L^{\phi}_{\infty}
\bigr)\nonumber\\
&&\qquad=\frac{\int
u(y,w)(\int
u(w,z)h_{x}(z) \,d\nu(z)) \,d\phi(w)} {h_{x}(y)}
\\
&&\qquad\quad{} + {\int u(y,w)(\int u(w,z)h_{x}(z) \,d\phi(z)) \,d\nu(w)
\over h_{x}(y)}.
\nonumber
\end{eqnarray}
By assumption (i) these are finite, and since they are excessive in
$y$ it follows that $H^{\nu}_t$ is right continuous in $t$. Hence, it
follows from (\ref{martj}) that
$A_t^\nu, t\geq0$, is also right continuous.
Therefore,
%
\begin{equation}
\mathop{\sup_{d(\nu,\nu')\le\delta}}_ {\nu,\nu'\in F}  A^\nu _t-A^{\nu'}_t
=\mathop{\sup_{d(\nu,z)\le\delta}}_ {\nu,\nu'\in F} \bigl|A^\nu _t-A^{\nu'}_t\bigr|
\end{equation}
is a right continuous, nonnegative submartingale and, therefore, using
(\ref{martj}), we see that
%
\begin{eqnarray} \label{12.0}
 I_{1,x}&=& P_{\phi}^{x/h_{x}} \Bigl( \sup
_{t\geq0} \mathop{\sup_{d(\nu
,\nu
')\le\delta}}_ {\nu,\nu'\in F} A^\nu_t-A^{\nu'}_t
\ge\ep\Bigr)
\nonumber
\\[-8pt]
\\[-8pt]
\nonumber
& \le&\frac{1}{\ep}E_{\phi}^{x/h_{x}}\Bigl( \mathop{\sup
_{d(\nu,\nu
')\le\delta}}_
{\nu,\nu'\in D} \bigl| L^{\nu}_{\infty}-L^{\nu'}_{\infty}\bigr|
\Bigr).
\end{eqnarray}

Using (\ref{nit.1+}) and then (\ref{de.1}),
%
\begin{eqnarray}\label{12.2aq}
&& Q^{\rho,h}_{\phi} \Bigl( \sup_{t\geq0}\mathop{\sup
_{d(\nu,\nu')\le
\delta}}_
{\nu,\nu'\in F}  L^{\nu}_{t}-L^{\nu'}_{t}
\ge3\ep\Bigr)
\nonumber\\
&&\qquad =\int Q_{\phi}^{x,h_{x}} \Bigl( \sup_{t\geq0}
\mathop{\sup_{d(\nu
,\nu')\le\delta}}_ {\nu,\nu'\in F} L^{\nu}_{t}-L^{\nu'}_{t}
\ge3\ep\Bigr) \,d\rho(x)
\\
&&\qquad =\int P_{\phi}^{x/h_{x}} \Bigl( \sup_{t\geq0}
\mathop{\sup_{d(\nu
,\nu')\le\delta}}_ {\nu,\nu'\in F} L^{\nu}_{t}-L^{\nu'}_{t}
\ge3\ep\Bigr) Q^{x,h_{x}}\bigl(L^{\phi}_{\infty}\bigr) \,d\rho
(x),
\nonumber
\end{eqnarray}
so that by (\ref{12.2a}) and (\ref{12.0})
%
\begin{eqnarray}\label{12.2ar}
&&Q^{\rho,h}_{\phi} \Bigl( \sup_{t\geq0}
\mathop{\sup_{d(\nu
,\nu
')\le\delta}}_ {\nu,\nu'\in F} L^{\nu}_{t}-L^{\nu'}_{t}
\ge3\ep\Bigr)
\nonumber
\\
&&\qquad\leq\frac{1}{\ep}\int E_{\phi}^{x/h_{x}}\Bigl(\mathop{ \sup
_{d(\nu
,\nu
')\le\delta}}_ {\nu,\nu'\in D} \bigl| L^{\nu}_{\infty}-L^{\nu
'}_{\infty}\bigr|
\Bigr) Q^{x,h_{x}}\bigl(L^{\phi}_{\infty}\bigr) \,d\rho(x)
\nonumber
\\[-8pt]
\\[-8pt]
\nonumber
&&\qquad\quad{} +\int P_{\phi}^{x/h_{x}} \Bigl(\sup_{t\geq0}
\mathop{\sup_{d(\nu
,\nu')\le\delta}}_ {\nu,\nu'\in F} H^\nu_t-H^{\nu'}_t
\ge2\ep\Bigr) Q^{x,h_{x}}\bigl(L^{\phi}_{\infty}\bigr) \,d\rho
(x)
\\
&&\qquad:=I_{\delta}+\mathit{II}_{\delta}.\nonumber
\end{eqnarray}
Using (\ref{de.1}) and then (\ref{nit.1+}), we see that
%
\begin{equation}
I_{\delta}=\frac{1}{\ep} Q^{\rho,h}_{\phi}\Bigl( \mathop{\sup
_{d(\nu
,\nu')\le
\delta}}_ {\nu,\nu'\in D}  \bigl| L^{\nu}_{\infty}-L^{\nu'}_{\infty
}\bigr|
\Bigr). \label{5.27}
\end{equation}
It follows from (\ref{low.2}) and
assumption ii
that for any $\ep'>0$, we can choose a
$\delta>0$, for which (\ref{5.27}) is less that $\ep'$.

We show below that
%
\begin{equation}
\lim_{\delta\rar0}P_{\phi}^{x/h_{x}} \Bigl(\sup
_{t\geq0} \mathop{\sup_{d(\nu,\nu')\le\delta}}_ {\nu,\nu'\in D} H^\nu_t-H^{\nu'}_t
\ge2\ep\Bigr)Q^{x,h_{x}}\bigl(L^{\phi}_{\infty}\bigr)=0,
\label{martk}
\end{equation}
uniformly in $x$.
Considering (\ref{12.2ar}), the proof is completed by taking
$F\uparrow D$.

To prove (\ref{martk}), we use (\ref{12.2as}) to write
%
{\fontsize{10.5}{12.5}{\selectfont
\begin{eqnarray}\label{martn}
&& P_{\phi}^{x/h_{x}} \Bigl(\sup_{t\geq0} \mathop{\sup
_{d(\nu,\nu')\le
\delta
}}_ {\nu,\nu'\in D} H^\nu_t-H^{\nu'}_t
\ge2\ep\Bigr) \nonumber
\\
&& \qquad\leq P_{\phi}^{x/h_{x}} \biggl(\sup_{t\geq0}
\mathop{\sup_{d(\nu
,\nu
')\le\delta}}_ {\nu,\nu'\in D} {h_{x}(X_{t}) L^{\phi}_{t} (E^{X_{t}/h_{x}}(L^{\nu}_{\infty
})-E^{X_{t}/h_{x}}(L^{\nu' }_{\infty})) \over
h_{x}(X_{t})E^{x/h_{x}}(L^{\phi}_{\infty}\mid\mathcal F^{0}_t)}\ge \ep\biggr)
\nonumber\\[-8pt]\\[-8pt]\nonumber
&&\qquad\quad{} + P_{\phi}^{x/h_{x}} \biggl(\sup_{t\geq0} \mathop{\sup
_{d(\nu,\nu
')\le\delta}}_ {\nu,\nu'\in D} {h_{x}(X_{t}) (E^{X_{t}/h_{x}}(L^{\nu}_{\infty}L^{\phi}_{\infty}
)-E^{X_{t}/h_{x}}(L^{\nu' }_{\infty}L^{\phi}_{\infty} ) ) \over
h_{x}(X_{t})E^{x/h_{x}}(L^{\phi}_{\infty}\mid\mathcal F^{0}_t)}\ge \ep\biggr)
\nonumber.
\end{eqnarray}}}
\hspace*{-2pt}Let
\[
\ga_{x}(\delta)=\sup_{y\in S} \mathop{\sup
_{d(\nu,\nu')\le
\delta}}_ {\nu
,\nu'\in D} h_{x}(y)\bigl|E^{y/h_{x}}
\bigl(L^{\nu}_{\infty}\bigr)-E^{y/h_{x}}\bigl(L^{\nu'
}_{\infty}
\bigr)\bigr|
\]
and
%
\begin{equation}
\bar\ga_{x}(\delta)=\sup_{y\in S} \mathop{\sup
_{d(\nu,\nu
')\le\delta}}_
{\nu,\nu'\in D} h_{x}(y)\bigl |E^{y/h_{x}}
\bigl(L^{\nu}_{\infty}L^{\phi}_{\infty}
\bigr)-E^{y/h_{x}}\bigl(L^{\nu' }_{\infty
}L^{\phi}_{\infty}
\bigr)\bigr |\label{12.2bv}.
\end{equation}
Then the first line of (\ref{martn}) is less than or equal to
%
\begin{eqnarray}\label{5.31}
&& P_{\phi}^{x/h_{x}} \biggl(\sup_{t\geq0}
{L^{\phi}_{t} \ga_{x}(\delta) \over h_{x}(X_{t})
E^{x/h_{x}}(L^{\phi}_{\infty
}\mid\mathcal F^{0}_t)}\ge\ep\biggr)
\nonumber
\\[-8pt]
\\[-8pt]
\nonumber
&&\qquad{} + P_{\phi}^{x/h_{x}} \biggl(\sup_{t\geq0}
{\bar\ga_{x}(\delta) \over h_{x}(X_{t}) E^{x/h_{x}}(L^{\phi
}_{\infty}\mid
\mathcal F^{0}_t)}\ge\ep\biggr)
.
\end{eqnarray}

It follows from (\ref{mema}), (\ref{memb}) and assumption (i) that
%
\begin{equation}
\lim_{\delta\rar0}\ga_{x}(\delta)=0 \quad\mbox{and}\quad \lim
_{\delta
\rar0}\bar\ga _{x}(\delta)=0\label{memc}
\end{equation}
uniformly in $x\in S$. Consequently, bounding $L^{\phi}_{t}$ by
$E^{x}(L^{\phi}_{\infty}\mid\mathcal F^{0}_t)$ in the first line of
(\ref
{5.31}), we see that (\ref{martk}) follows from the next lemma.
\end{pf}

\begin{lemma}\label{lem-supmart}
Let $M_{t}$ be a nonnegative right continuous $P^{x}$ martingale. Then
%
\begin{equation}
\frac{M_{t}}{h_{x}(X_{t})E^{x/h_{x}}(L^{\phi}_{\infty}\mid\mathcal
F^{0}_t)}, \qquad t\geq0\label{5.34}
\end{equation}
is a right continuous nonnegative
supermartingale with respect to $P_{\phi}^{x/h_{x}} $,
and
%
\begin{equation}
P_{\phi}^{x/h_{x}}\biggl(\sup_{t\geq0}
{M_{t} \over
h_{x}(X_{t})E^{x/h_{x}}(L^{\phi}_{\infty}\mid\mathcal F^{0}_t) }\ge \ep \biggr)\le \frac{1}{\ep}
{P^{x}(M_{0}) \over Q^{x,h_{x}}(L^{\phi}_{\infty
})}.\label{12.max}
\end{equation}
\end{lemma}

\begin{pf} For any $t>s\geq0 $ and any $F_s\in\mathcal{F}^{0}_{s}$,
we have
%
\begin{eqnarray}\label{j6.10a}
J&:=&P_{\phi}^{x/h_{x}}\biggl( F_s
\frac
{M_{t}}{h_{x}(X_{t})E^{x/h_{x}}(L^{\phi}_{\infty}\mid\mathcal F^{0}_t)
}\biggr)\nonumber
\\
&=&{1 \over E^{x/h_{x}}(L^{\phi}_{\infty})}P^{x/h_{x}}\biggl( L_{\infty
}^{\phi
}F_s
\frac{M_{t}}{h_{x}(X_{t})E^{x/h_{x}}(L^{\phi}_{\infty}\mid
\mathcal
F^{0}_t) }\biggr)
\\
&=& {1 \over E^{x/h_{x}}(L^{\phi}_{\infty})}P^{x/h_{x}}\biggl( F_s
\frac
{M_{t}}{h_{x}(X_{t})}\biggr).
\nonumber
\end{eqnarray}
Note that for all functions $f$ on $S$, $f(\De)=0$. Therefore, using
(\ref{nit4.4a}) and (\ref{low.3})
%
\begin{equation}\qquad
P^{x/h_{x}}\biggl( F_s\frac{M_{t}}{h_{x}(X_{t})}\biggr)=P^{x/h_{x}}
\biggl( 1_{\{\ze>t\}
}F_s\frac{M_{t}}{h_{x}(X_{t})}\biggr)=
{P^{x}( 1_{\{\ze>t\}}F_sM_{t}) \over
h_{x}(x) }.\label{j6.10aq}
\end{equation}
Consequently,
\begin{eqnarray*}
J&=&{1 \over h_{x}(x) E^{x/h_{x}}(L^{\phi}_{\infty})}P^{x}( 1_{\{\ze
>t\}
}F_sM_{t})
\\
&\leq&{1 \over h_{x}(x) E^{x/h_{x}}(L^{\phi}_{\infty})}P^{x}( 1_{\{
\ze
>s\}}F_sM_{t})
\\
&=&{1 \over h_{x}(x) E^{x/h_{x}}(L^{\phi}_{\infty})} P^{x}(1_{\{\ze
>s\}
}F_sM_{s})
.
\end{eqnarray*}
Considering (\ref{j6.10aq}) and (\ref{j6.10a}) with $t$ replaced by
$s$, we see that the last line above is equal to
%
\begin{equation}
P_{\phi}^{x/h_{x}}\biggl(F_s\frac{M_{s}}{h_{x}(X_{s})E^{x/h_{x}}(L^{\phi
}_{\infty}\mid\mathcal F^{0}_s)}
\biggr).
\end{equation}
This shows that (\ref{5.34}) is a nonnegative
supermartingale with respect to $ P_{\phi}^{x/h_{x}}$. That it is right
continuous follows from
%
\begin{equation}
E^{x/h_{x}}\bigl(L^{\phi}_{\infty}\mid\mathcal
F^{0}_t\bigr) = L^{\phi
}_{t}+E^{X_{t}/h_{x}}
\bigl(L^{\phi}_{\infty} \bigr)
\end{equation}
and the sentence following (\ref{memb}). This and the fact that
$h_{x}(x) E^{x/h_{x}}(L^{\phi}_{\infty})= Q^{x,h_{x}}(L^{\phi
}_{\infty})$
gives (\ref{12.max}).
\end{pf}

We can now give our most general result about the joint continuity of
the continuous additive functionals.

\begin{theorem}\label{theo-pbound} Assume that conditions \textup{(i)} and \textup{(ii)}
in Lemma~\ref
{nit1} are satisfied for some $ \phi\in\mathcal{R}_{\|\cdot\|
_{u^2,\infty
}}^{+}$ with support $\phi=S$, and some $\rho\in\mathcal{R}_{\|
\cdot\|
_0}^{+}$ of the form $\rho(dx)=f(x) m(dx)$ with $f>0$. Then there exists
a version of $ \{L^{\nu}_{t}, (t,\nu)\in R^{1}_{+}\times\VV
\} $
that is continuous on $(0,\ze)\times\VV$, $P^{x}$ almost surely for
all $x\in S$, and is continuous on $[0,\ze)\times\VV$, $P^{x}$ almost
surely for $m(dx)$ a.e. $x\in S$. (Continuity on $\VV$ is with respect
to the metric $d$ introduced in the statement of Lemma~\ref{nit1}.)
\end{theorem}

\begin{pf} The first step in this proof is to show that
$ \{L^\nu_t, (t,\nu)\in R^{+}\times D
\}$ is locally uniformly continuous almost surely with respect to $
Q_{\phi}^{\rho,h}$. This can be proved by mimicking the proof in
\cite{sip}, Theorem~6.1, that (6.10) implies (6.12). (This theorem is given
for a different family of continuous additive functionals with
different conditions on the
potential density of the associated Markov process, nevertheless it is
not difficult to see that a straightforward adaptation of the proof
works in the case we are considering.)

Let
%
\begin{equation}
\tilde{\Om}_{1}= \bigl\{ \omega| L^\nu_t(
\omega) \mbox{ is locally uniformly continuous on } R^{+}\times D \bigr
\}.\label{om1}
\end{equation}
We have that
%
\begin{equation}
Q_{\phi}^{\rho,h}\bigl(\tilde{\Om}_{1}^{c}
\bigr)=\int Q_{\phi
}^{x,h_{x}}\bigl(\tilde {\Om}_{1}^{c}
\bigr) \,d\rho(x) =0.\label{2.5}
\end{equation}
Using the fact that $L^{\phi}_{\infty}>0$ and $\rho(dx)=f(x) m(dx)$ with
$f>0$, we see from~(\ref{2.5}) that
%
\begin{equation}
Q^{x,h_{x}}\bigl(\tilde{\Om}_{1}^{c}\bigr) =0\qquad\mbox{for $m(dx)$ a.e. $x\in S$.}\label{2.5q}
\end{equation}
Set
%
\begin{equation}\quad
\tilde{\Om}_{2}= \bigl\{ \omega| L^\nu_t(
\omega) \mbox{ is locally uniformly continuous on } [0,\ze)\times D \bigr\}
.\label{om2}
\end{equation}
We see from (\ref{2.5q}) and (\ref{nit4.4a}) that
%
\begin{equation}
P^{x}\bigl(\tilde{\Om}_{2}^{c}\bigr) =0\qquad \mbox{for $m(dx)$ a.e. $x\in S$.}\label{2.5r}
\end{equation}
Because the Markov process has transition densities, we see that for
any $x\in S$ and $\ep>0$
%
\begin{equation}
P^{x}\bigl(\tilde{\Om}_{2}^{c}\circ
\th_{\ep}\bigr)=E^{x}\bigl(P^{X_{\ep}}\bigl(\tilde {\Om
}_{2}^{c}\bigr) \bigr)=\int p_{\ep}(x,y)
P^{y}\bigl(\tilde{\Om}_{2}^{c}\bigr) \,dm(y)=0.
\label{2.5s}
\end{equation}
Consequently, for
%
\begin{equation}
\qquad \tilde{\Om}_{3}:= \bigl\{ \omega| L^\nu_t(
\omega) \mbox{ is locally uniformly continuous on } (0,\ze)\times D \bigr\}
,\label{om3}
\end{equation}
we have
%
\begin{equation}
P^{x}\bigl(\tilde{\Om}_{3}^{c}\bigr) =0\qquad
\mbox{for all $x\in S$.}\label{2.5t}
\end{equation}

For $\omega\in\tilde{\Om}_{3}^{c}$ we set $\tilde L^\nu_t
(\omega)\equiv
0$. For $\omega\in\tilde{\Om}_{3}$ we define $ \{\tilde L^\nu_t
(\omega),
(t,\nu)\in(0,\ze)\times\VV
\}$ as the continuous extension of
$ \{L^\nu_t (\omega), (t,\nu)\in(0,\ze)\times D
\}$, and then set
%
\begin{equation}
\tilde L^\nu_0(\omega)=\mathop{\liminf_{s\downarrow0}}_{s\ \mathrm{rational}}
\tilde L^\nu_s(\omega)\label{5.42a}
\end{equation}
and
%
\begin{equation}
\tilde L^\nu_t(\omega)= \mathop{\liminf_{s\uparrow\zeta(\omega)}}_{s\ \mathrm{
rational}}
\tilde L^\nu_s(\omega)\qquad\mbox{for all }t\geq\ze
.\label{5.42}
\end{equation}

Since $ L^\nu_t (\omega)$ is increasing in $t$ for $\nu\in D$, the same
is true for $ \{\tilde L^\nu_t (\omega), (t,\nu)\in(0,\ze)\times
\VV
\}$. Therefore the lim infs in (\ref{5.42a}) and (\ref{5.42}) are
actually limits.
Since we can assume that the
$L^\nu_t$ are perfect continuous additive functionals for all
$\nu\in D$, we immediately see that the same is true for
$\tilde L^\nu_t$ for each $\nu\in\VV$, except that one problem
remains. We need $\tilde L^\nu_0=0$, but it is not clear from (\ref
{5.42a}) that this is the case.

We show that
$\tilde{L}^{\nu}_{t}$ is a version of $L^{\nu}_{t}$, which implies
that $\tilde L^\nu_0=0$.
Pick some $\nu'$ not in $D$ and set $D'=D\cup\{\nu'\}$. Then by the
argument above, but with $D$ replaced by $D'$, we get that
$L^\nu_t(\omega) \mbox{ is locally uniformly continuous on}
(0,\ze)\times D'$, almost surely. Thus, $L^{\nu'}_t=\widetilde
L^{\nu'}_t$ on
$(0,\ze)$ a.s., which is enough to show that $\{\widetilde L^{\nu
'}_t, t\geq
0\}$ is a version of $\{L^{\nu'}_t, t\geq0\}$.

Thus, we see that there exists
a version of $ \{L^{\nu}_{t}, (t,\nu)\in R^{1}_{+}\times\VV
\} $
that is continuous on $(0,\ze)\times\VV$, $P^{x}$ almost surely for
all $x\in S$.
To see that this version is continuous on $[0,\ze)\times\VV$, $P^{x}$
almost surely for $m(dx)$ a.e. $x\in S$, it suffices to note that for
each $\omega\in\Om_{2}$, $\tilde L^\nu_t (\omega)$ is continuous
on $[0,\ze
)\times\VV$, and then use (\ref{2.5r}).
\end{pf}

We now take $S=R^n$. Let $T_{a}$ denote the bijection on the space of
measures defined by the translation $T_{a}(\nu)=\nu_{a}$; see (\ref
{1.19}). We say that a set $\VV$ of measures on $R^{n}$ is translation
invariant if it is invariant under $T_{a}$ for each $a\in R^{n}$
and say that a topology $\OO$ on such a set $\VV$
is homogeneous if $T_{a}$ is an isomorphism for each $a\in R^{n}$.

\begin{theorem}\label{t1.3} Let $X$ be an exponentially killed L\'evy
process in
$R^{n}$ and $\VV\subseteq\RR_{\|\cdot\|}^{+} $ be a translation
invariant set of measures on $R^{n}$. Assume:
\begin{longlist}[{(ii)}]
\item[{(i)}]
that there is a
homogeneous topology $\OO$ for $\VV$ under which $\VV$ is a separable
locally compact metric space with metric $d$, and
\item[{(ii)}]that conditions \textup{(i)} and \textup{(ii)} in Lemma~\ref{nit1}
are satisfied for some $ \phi\in\mathcal{R}_{\|\cdot\|_{u^2,\infty
}}^{+}$ with support $\phi=S$, and some $\rho\in\mathcal{R}_{\|
\cdot\|
_0}^{+}$ of the form $\rho(dx)=f(x) m(dx)$ with $f>0$.
\end{longlist}

Then there exists
a version of $ \{ L^{\nu}_{t}, (t,\nu)\in R^1_{+}\times
\VV \} $ that is continuous $P^{x}$ almost surely for all $x\in S$.
\end{theorem}

\begin{pf} Using the fact that $X$ is an exponentially killed process, it
follows easily from the proof of Theorem~\ref{theo-pbound} and \cite{MR96}, page~1149, that we can replace $\ze$ by $\infty$ in the conclusions of
Theorem~\ref{theo-pbound}. Hence, there exists a version of $ \{
L^{\nu}_{t}, (t,\nu)\in R^1_{+}\times
\VV \} $ that is continuous $P^{x}$ almost surely for a.e. $x\in S$.
By translation invariance, this holds for all $x\in S$.
\end{pf}

By Corollary~\ref{cor-4.1}, we can replace condition (ii) of Lemma~\ref
{nit1} by (\ref{6p}). This is used to obtain the next corollary that
allows us to replace condition (ii) in Lemma~\ref{nit1} by a more
concrete condition that follows from Theorem~\ref{theo-contprop}.

\begin{corollary}
\label{nit1af}Let $X =
(\Om, X_t,P^x
)$ be a Borel right process in $S$ with strictly positive
$0$-potential densities
$u(x,y)$, and let
$\VV$ be a separable locally compact subset of $ \mathcal{R}_{\|\cdot
\|}^{+}$.
Assume that there exists
a probability measure $\si$ on $\VV$ such that $J_{\VV, \|\cdot\|
,\si
}(D)<\infty$, where $D$ is the diameter of $\VV$ with respect to $ \|
\cdot
\|$, and
%
\begin{equation}
\lim_{\delta\to0} J_{\VV, \|\cdot\|,\si}(\delta)=0. \label{2.10v}
\end{equation}
Then condition \textup{(ii)} of Lemma~\ref{nit1} holds.
\end{corollary}

\section{Continuous additive functionals of L\'evy processes}\label{sec-levy}

The main purpose of this section is to prove Theorem~\ref{theo-1.1}.
We begin with two lemmas which follow easily from results in \cite
{sip}. Because the notation in \cite{sip} is different from the
notation used in Theorem~\ref{theo-1.1} it is useful to be more
explicit about the relationship between a L\'evy process killed at the
end of an independent exponential time and the L\'evy process itself,
that is, the unkilled process. Let
$Y=\{Y_{t},t\in R^{+}\}$ be a L\'evy process in $R^d$ with
characteristic exponent $\bar\ka$. Let $X=\{X_{t},t\in R^{+}\}$ be the
process $Y$, killed at the end of an independent exponential time with
mean $1/\bb$. Let $\ka$ and $u$ denote the characteristic exponent of~$X$ and the
potential density of~$X$.
Then
%
\begin{equation}
\ka(\xi)=\bb+\bar\ka(\xi) \label{link1}
\end{equation}
and
%
\begin{equation}
\widehat u (\xi)={1 \over\ka(\xi)}={1 \over\bb+\bar\ka(\xi
)}.\label{link2}
\end{equation}

\begin{lemma}\label{lem-7.1} Let $X$ be a L\'evy process in $R^{d}$
that is killed at the end of an independent exponential time, with
characteristic exponent $\ka$ and potential density~$u$, and suppose that
%
\begin{equation}
\frac{1}{|\ka(\xi)|^{2}}\le C\frac{\ga(\xi)}{|\xi|^{d}},\label{7.6w}
\end{equation}
where $\ga=|\widehat {u}|\ast|\widehat { u}|$. Then
%
\begin{equation}
\int\bigl|\hat\nu(s)\bigr|\bigl|\widehat{u}(s)\bigr| \,ds\le C \int_{1}^{\infty}
{(\int_{|\xi
|\geq x}\ga(\xi)|\hat{\nu}(\xi)|^{2} \,d\xi) ^{1/2}\over
x(\log2x)^{1/2}} \,dx.\label{7.6}
\end{equation}
\end{lemma}

\begin{pf} We follow the proof of \cite{sip}, Lemma~5.2, with the $\ga
$ of
this theorem and $|\ka(\xi)|$ replacing the $\ga$ and $(1+\psi(\xi))$
in \cite{sip}, Lemma~5.2. It is easy to see that the proof of
\cite{sip}, Lemma~5.2, goes through with these changes to prove this
lemma.
\end{pf}

\begin{remark}\label{rem-7.1}
Since
%
\begin{equation}
\sup_{y}\bigl|U\nu(y)\bigr|\leq C\int\bigl|\hat\nu(s)\bigl|\bigl|\widehat{u}(s)\bigr|
\,ds,\label{bbnd7}
\end{equation}
it follows from (\ref{7.6}) and \cite{BG}, page~285, that $\nu$ charges no
polar set. It is a conjecture of Getoor that essentially all L\'evy
processes in $R^{d}$ satisfy Hunt's hypothesis (H) which is that all
semipolar sets are polar. This has been proved in many cases. See, for
example, \cite{Rao} and \cite{GR}.
In these cases, the condition in Theorem~\ref{theo-1.1}, that
$\nu\in\mathcal{R}^{+}(X)$, is superfluous.
\end{remark}

\begin{remark}\label{rem-7.2} The function $\ga(\xi)$ plays a
critical role in Theorem~\ref{theo-1.1}. We note that
%
\begin{equation}
\sup_{\xi\in R^{d}}\ga(\xi)<C\|u\|_{2}^{2}
\end{equation}
for some absolute constant $C$.
\end{remark}

The next lemma is a generalization of \cite{sip}, Lemma~5.3.

\begin{lemma} \label{lem-7.2}If
%
\begin{equation}
C _{1}\tau\bigl(|\xi|\bigr) \le\bigl|\ka(\xi)\bigr|\le C _{2}\tau\bigl(|\xi|\bigr)\qquad
\forall\xi \in R^{d} \label{7.8}
\end{equation}
and $\tau(|\xi|) $ is regularly varying at infinity, then (\ref
{7.6w}) holds.
\end{lemma}

\begin{pf} By the assumption of regular variation, for $|\xi|$
sufficiently large,
%
\begin{eqnarray}\label{7.5q}
\ga(\xi)&\ge& \int_{|\eta|\ge2|\xi|}\frac{d\eta}{|\ka(\xi
-\eta)||\ka(
\eta)|}\nonumber
\\[-2pt]
&\ge& \int_{|\eta|\ge2|\xi|}\frac{d\eta}{ \tau( |\eta-\xi|)
\tau(
|\eta|) }
\nonumber
\\[-9pt]
\\[-9pt]
\nonumber
&\ge& \int_{|\eta|\ge2|\xi|}\frac{d\eta}{ \tau^{2}( |\eta|)
}
\nonumber
\\[-2pt]
& \ge&
\nonumber
C\frac{|\xi|^{d}}{\tau^{2}(|\xi|)},
\end{eqnarray}
which gives (\ref{7.6w}).
(Since this is a lower bound, it holds even if the integral on the
third line is infinite.)
It is clear that
the constant in (\ref{7.5q}) can be adjusted to hold for all $\xi\in
R^{d}$.
\end{pf}

The following lemma provides a key estimate in the proof of Theorem~\ref{theo-1.1}.

\begin{lemma} \label{lem-7.3a} Under the hypotheses of Lemma~\ref{lem-7.2},
%
\begin{equation}
\int\bigl|\widehat{u}(\la_{1})\bigr| ^{2}\bigl|\widehat{u}(\xi-
\la_{1})\bigr| \,d\la_{1} \le C \bigl|\widehat{u}(\xi)\bigr| \|u
\|_{2}^{2}.\label{7.12}
\end{equation}\vadjust{\goodbreak}
\end{lemma}

\begin{pf}
Using (\ref{7.8}), we can treat $u$ as though $|\hat u(|\xi|)|$ is
regularly varying at infinity. Consequently,
%
\begin{eqnarray}\label{simpy}
&& \int\bigl|\widehat{u}(\la_{1})\bigr| ^{2}\bigl|\widehat{u}(\xi-
\la_{1})\bigr| \,d\la _{1}\nonumber
\\
&&\qquad \le\int_{|\la_{1}\le|\xi|/2|}\bigl |\widehat{u}(\la_{1})\bigr|
^{2}\bigl|\widehat {u}(\xi-\la_{1})\bigr| \,d\la_{1}+\int
_{|\la_{1}\ge|\xi|/2|} \bigl|\widehat {u}(\la _{1})\bigr| ^{2}\bigl|
\widehat{u}(\xi-\la_{1})\bigr| \,d\la_{1}
\nonumber
\\[-8pt]
\\[-8pt]
\nonumber
&&\qquad
\le C \bigl|\widehat{u}(\xi)\bigr|\biggl(\int_{|\la_{1}\le|\xi|/2|} \bigl|
\widehat {u}(\la_{1})\bigr| ^{2} \,d\la_{1}+ \int
_{|\la_{1}\ge|\xi|/2|} \bigl|\widehat {u}(\la_{1})\bigr| \bigl|\widehat{u}(\xi-
\la_{1})\bigr| \,d\la_{1}\biggr)
\\
&&\qquad
\le C\bigl |\widehat{u}(\xi)\bigr|\bigl( \| {u} \|^{2}_{2}+
\ga(\xi)\bigr),\nonumber
\end{eqnarray}
which implies (\ref{7.12}).
\end{pf}

\begin{pf*}{Proof of Theorem~\ref{theo-1.1}} This theorem is an
immediate consequence of Theorem~\ref{t1.3}. We begin by showing that
Theorem~\ref{t1.3}(ii) holds.
We take $\phi(dx)=\rho(dx)=e^{-|x|^{2}/2} \,dx$ and we set
$h(y,x)=h(x-y)=u_{1}(x-y)$ where $u_{1}(y)=\int_{1}^{\infty}
p_{t}(y) \,dt$.
We have
%
\begin{equation}
\bigl|\hat h (\la)\bigr|= \bigl|\widehat {u} (\la)\bigr| e^{- \operatorname{Re} \kappa(\la)}.\label{hhat}
\end{equation}

To show that condition (i) of Lemma~\ref{nit1} holds we show that
%
\begin{equation}
\sup_{x,y} \biggl|\int u (y,z) h(z,x) \,d\nu(z) \biggr|\le C \int\bigl|\hat \nu
(s)\bigr|\bigl|\widehat {u}(s)\bigr| \,ds,\label{m6.11q}
\end{equation}
and
%
\begin{equation}
\qquad\sup_{x,y}\biggl | \int u(y,z) \biggl(\int u(z,w)h(w,x) \,d\phi(w)
\biggr) \,d\nu (z)\biggr |\le C \int\bigl|\hat\nu(s)\bigr| \bigl|\widehat {u}(s)\bigr| \,ds.\label{m6.11aq}
\end{equation}
When (\ref{b101}) holds, it follows from (\ref{7.6}) that the
right-hand side is finite.
Therefore,
replacing $\nu$ in (\ref{m6.11q}) and (\ref{m6.11aq}) by $\nu
_{r}-\nu
_{r'}$, so that
$|\hat\nu(s)|$ is replaced by $|e^{ir\cdot s}-e^{ir'\cdot s} | |\hat
\nu(s)|$, we see that condition (i) of Lemma~\ref{nit1} follows from
(\ref{7.6}) and the dominated convergence theorem.

To obtain (\ref{m6.11q}), we write
%
\begin{eqnarray}
&& \int u(y,z) h(z,x) \,d\nu(z)
\nonumber\\
&&\qquad=
\int u(z-y) h(x-z) \,d\nu(z)
\nonumber
\\[-8pt]
\\[-8pt]
\nonumber
&&\qquad=
\int e^{i(z-y)\la_{1}}\widehat{u}(\la_{1})
e^{i(x-z)\la
_{2}}\hat h(\la_{2}) \,d\la_{1} \,d
\la_{2} \,d\nu(z)
\\
&&\qquad=
\int\hat\nu(\la_{1}-\la_{2})e^{-iy \la_{1}}
\widehat {u}(\la _{1}) e^{ix \la_{2}}\hat h(\la_{2}) \,d
\la_{1} \,d\la_{2}.\nonumber
\end{eqnarray}
Hence,
%
\begin{equation}\label{exc.2q}\qquad
\sup_{x,y} \biggl|\int u(y,z) h(z,x) \,d\nu(z) \biggr|  \leq
\int\bigl|\hat\nu(s)\bigr|\biggl(\int\bigl|\widehat{u}(s+\la_{2}) \bigr|\bigl|\hat h(
\la_{2})\bigr| \,d\la_{2}\biggr) \,ds.
\end{equation}
We complete the proof of (\ref{m6.11q}) by showing that
%
\begin{equation}
\int\bigl|\widehat{u}(s+\la)\bigr |\bigl|\hat h(\la)\bigr| \,d\la\le C\bigl |\widehat{u}(s)\bigr |.
\label{7.11}
\end{equation}
We have
%
\begin{equation}
\int\bigl|\widehat{u}(s+\la) \bigr|\bigl|\hat h(\la)\bigr| \,d\la=C \int\frac{e^{-\operatorname{Re}
\ka
(\la)}}{| \ka(s+\la)|| \ka( \la)|} \,d
\la.\label{7.17}
\end{equation}
Using the same inequalities used in the proof of Lemma~\ref{lem-7.3a},
we see that
%
\begin{equation}
\int_{|\la|\le|s|/2}\frac{e^{- \operatorname{Re} \ka(\la)}}{|\ka(s+\la)||\ka(
\la)|} \,d\la\le C\frac{1}{|\ka(s )|}
\int\frac{e^{- \operatorname{Re} \ka(\la)}}{ |\ka
( \la
)|} \,d\la\label{7.18}
\end{equation}
and
%
\begin{equation}
\int_{|\la|\ge|s|/2}\frac{e^{- \operatorname{Re} \ka(\la)}}{|\ka(s+\la)||\ka(
\la)|} \,d\la\le C\frac{1}{|\ka(s )|}
\int e^{- \operatorname{Re} \ka(\la)} \,d\la.\label{7.19}
\end{equation}
Using (\ref{expint}), (\ref{7.18}) and (\ref{7.19}) in (\ref{7.17}) and
then (\ref{1.19cor}), we get (\ref{7.11}).

In a similar manner to how we obtained (\ref{exc.2q}) by taking
Fourier transforms, we see that
%
\begin{eqnarray}\label{exc.13}
&& \sup_{x,y} \biggl| \int u(y,z) \biggl(\int u(z,w)h(w,x) \,d\phi(w)
\biggr) \,d\nu (z)\biggr |
\nonumber\\
&& \qquad\leq\int\biggl(\int\bigl|\hat\phi(\la_{1}-\la_{2})\bigr| \bigl|\hat h(
\la_{2})\bigr| \,d\la _{2}\biggr)\bigl| \widehat{u}(\la_{1})\bigr|
\bigl|\widehat{u}(\la_{3})\bigr| \bigl|\hat\nu\bigl((\la _{1}+
\la_{3})\bigr)\bigr| \,d\la_{1} \,d\la_{3}
\\
&& \qquad= \int\biggl(\int\bigl|\hat\phi(\la_{1}-\la_{2})\bigr| \bigl|\hat h(
\la_{2})\bigr| \,d\la _{2}\biggr)\bigl|\widehat{u}(\la_{1})\bigr|
\bigl|\widehat{u}(s-\la_{1})\bigr| \,d\la_{1} |\hat\nu (s)| \,ds.
\nonumber
\end{eqnarray}
Clearly, since $\phi=e^{-|x|^{2}/2} \,dx$, $|\hat\phi(\la)|\le
C|\widehat
{u}(\la)|$. Therefore, by (\ref{7.11})
%
\begin{equation}
\int\bigl|\hat\phi(\la_{1}-\la_{2})\bigr| \bigl|\hat h(
\la_{2})\bigr| \,d\la_{2}\le C\bigl|\widehat{u}(\la_{1})\bigr|.\label{7.23}
\end{equation}
Using this (\ref{exc.13}) and Lemma~\ref{lem-7.3a}, we get
(\ref{m6.11aq}).

We now show that condition (ii) of Lemma~\ref{nit1} holds. We have
already seen that $\nu\in\mathcal{R}^{+}$. Let
%
\begin{equation}
\|\nu\|_{\ga,2}: = \biggl(\int\bigl|\hat\nu(x )\bigr|^{2} \ga(x) \,dx
\biggr)^{1/2}.\label{bp.1q}
\end{equation}
It follows from \cite{LMR2}, Lemma~2.2, (see also \cite{LMR}, Theorem~6.1), that $\|\cdot\|_{\ga,2}$ is a proper norm for $u$, and it
follows from (\ref{b101}) that $\{\nu_{x}, x\in R^{d}\}\subseteq\RR
^{+}_{\|\cdot\|_{\ga.2}}$. To complete the proof of the continuity part
of Theorem~\ref{theo-1.1}, we need the following lemma which is proved below.

\begin{lemma}\label{lem-7.4} For any compact set $D\in R^{d}$,
%
\begin{equation}
\lim_{\delta\rar0}Q^{\rho}_{\phi} \Bigl(\mathop{\sup
_{|x-y|\le
\delta}}_ {\nu
,\nu'\in D} \bigl |L^{\nu_{x} }_{\infty} -L^{\nu_{y} }_{\infty}
\bigr|\Bigr)=0.\label{6pq}
\end{equation}
\end{lemma}

\begin{pf*}{Proof of Theorem~\ref{theo-1.1} continued} It
follows from Lemma~\ref{lem-7.4}
that condition (ii) of Lemma~\ref{nit1} holds with the metric $d$ being
the Euclidean metric on~$R^{d}$. Therefore, the conditions in Theorem~\ref{t1.3}(ii) hold and since $d$ is the Euclidean metric the condition
in Theorem~\ref{t1.3}(i)
also holds. The continuity portion of Theorem~\ref{theo-1.1} now
follows from Theorem~\ref{t1.3}.
\end{pf*}

\begin{pf*}{Proof of Lemma~\ref{lem-7.4}}This follows easily from
the proof of \cite{sip}, Theorem~1.6, with the $\ga$ of this theorem
replacing the $\ga$ in \cite{sip}, Theorem~1.6. The gist of the proof
of \cite{sip}, Theorem~1.6, is that (\ref{b101}) implies that for
compact sets $D$ of~$R^{d}$
%
\begin{equation}
\lim_{\delta\to0} J_{\overline \VV, \|\cdot\|_{\ga,2},\la
}(\delta)=0, \label{7.3}
\end{equation}
where $\overline \VV=\{\nu_{x},x\in D\}$ and $\la$ is Lebesgue
measure on
$R^{d}$. (See Section~\ref{sec-contperm} for notation.)

By Theorem~\ref{cor-jack}, we get that (\ref{6pq})
holds with $|x-y|$ replaced by $\|\nu_{x}-\nu_{y}\|_{\ga,2}$ and $
{x,y\in D}$ replaced by $ {\nu_{x},\nu_{y}\in\overline \VV}$. Since
%
\begin{equation}
\|\nu_{x}-\nu_{x+h}\|_{\ga,2}=C\biggl(\int
_{\xi\in R^{d}}\sin^{2}\frac
{\xi
h}{2}\ga(\xi)\bigl|\hat\nu(
\xi)\bigr|^{2} \,d\xi\biggr)^{1/2}\label{mjack}
\end{equation}
we see that $ \psi(\nu_{x})$ is continuous on $R^{d}$ and we get
(\ref
{6pq}) as stated.
\end{pf*}
\noqed\end{pf*}

Strengthening the hypotheses of Theorem~\ref{theo-1.1}, we get the
simple estimate of $\ga(\xi)$ in the next lemma.

\begin{lemma}\label{lem-7.3} Under the hypotheses of Lemma~\ref
{lem-7.2}, assume also that $\tau$ is regularly varying at infinity
with index greater than $ d/2 $ and less than $d$. Then
%
\begin{equation}
\ga(\xi)\le C \frac{|\xi|^{d}}{\tau^{ 2}(|\xi|) }\label{7.25}
\end{equation}
for all $|\xi|$ sufficiently large.
\end{lemma}

\begin{pf} This follows from \cite{LMR2}, Corollary~8.1.
\end{pf}

\begin{remark} We give some details on how the examples in Example~\ref{ex-1.1}, 1. and 2. are obtained.
\begin{longlist}[1.]
\item[1.] It is easy to see that (\ref{amp2}) follows from (\ref{b101}) and
Lemma~\ref{lem-7.3}.
\item[2.] In this case, the estimate in (\ref{7.25}) is not correct. To
find a bound for $\ga(\xi)$, we look at the proof of Lemma~\ref
{lem-7.3} with $d=2$ and $\tau$ as given in (\ref{1.25}). The bounds in
III remains the same but the bounds in I and II are now
%
\begin{equation}
C \frac{|\xi|^{2} \log|\xi|}{\tau^{2} (|\xi|) }
\end{equation}
for all $|\xi|$ sufficiently large. Given this, the rest of the
argument is essentially the same as in 1.
\end{longlist}
\end{remark}

We now take up the proof of the modulus of continuity assertion in
Theorem~\ref{theo-1.1}. We begin with a modulus of continuity result
for certain permanental processes, including the those considered in
Theorem~\ref{theo-1.1}.

\begin{theorem}\label{l5.1} Let $\{\psi(\nu),\nu\in\VV\}$ be a
permanental process with kernel $u$, where
$\VV=\{\nu_{x},x\in R^n\}$ is a family of measures such that
%
\begin{equation}
\| \nu_{x}-\nu_{y} \|\le\varrho\bigl(|x-y|\bigr),\label{7.28}
\end{equation}
where $\| \cdot\|$ is a proper norm on $\VV$ with respect to $u$, and
$\varrho$ is a strictly increasing function. Let
%
\begin{equation}
\omega(\delta)= \varrho(\delta)\log1/\delta+\int_{0}^{\delta}
\frac{ \varrho(u)}{u} \,du,\label{7.29}
\end{equation}
and assume that the integral is finite. Then for each $K>0$ there
exists a constant $C$ such that
%
\begin{equation}
\limsup_{\delta\to0 }\mathop{\sup_{ |x-y| \le\delta}}_{ x,y\in[-K,K]^{n}} \frac{\psi(\nu_{x})-\psi(\nu_{y}) }{\omega(\delta) }
\le C \qquad\mbox{a.s.}\label{15.1az}
\end{equation}

In particular, if $\varrho$ is a regularly varying function at zero
with index greater than zero,
we can take
%
\begin{equation}
\omega(\delta)= \varrho(\delta)\log1/\delta.\label{7.31}
\end{equation}
\end{theorem}

\begin{pf} This is proved in \cite{book}, Section~7.2, in a slightly
different setting. For it to hold in our setting, just change $(\log
1/u)^{1/2}$ in \cite{book}, (7.90), to $ \log1/u $ and continue the
proof with this change. This takes into account the fact that in (\ref
{tau}) we have a $\log$ rather than $(\log)^{1/2}$, which is what we
have when dealing with Gaussian processes.
\end{pf}

\begin{example}\label{ex-7.1} We consider Theorem~\ref{l5.1} in
the case where $u(x,y)=u(y-x)$ and the proper norm is $\|\cdot\|_{\ga
,2}$. By (\ref{bp.1q})
%
\begin{eqnarray}\label{7.32}
\|\nu_{x}-\nu_{y}\|_{\ga,2}&\le&C\biggl(\int\bigl|\hat
\nu_{x}(\la)-\nu _{y}(\la )\bigr|^{2} \ga(\la) \,d\la
\biggr)^{1/2}\nonumber
\\
&\le&C\biggl(\int\sin^{2}\frac{ (x-y)\la}2\bigl| \hat\nu(
\la)\bigr|^{2} \ga (\la) \,d\la\biggr)^{1/2}
\\
&\le&C \varphi\bigl(|x-y|\bigr),
\nonumber
\end{eqnarray}
where $\varphi$ is given in (\ref{1.23}). Note that if (\ref{b101})
holds then $\int(\varphi(u)/u) \,du<\infty$. Therefore, if (\ref{b101})
holds, the results in (\ref{15.1az})--(\ref{7.31}) hold with $\varrho$
replaced by~$\varphi$.
\end{example}

\begin{pf*}{Proof of Theorem~\ref{theo-1.1}, modulus of continuity}
This follows from Theorem~\ref{l5.1}, Example~\ref{ex-7.1} and the
second isomorphism theorem, Theorem~\ref{theo-ljr}, as in the proof of
a similar result in \cite{MR96}, Section~7. Note that the requirement
that $U^{1}\mu<\infty$ in \cite{MR96}, Theorem~2.2, follows from
(\ref
{b101}), (\ref{7.6}) and (\ref{bbnd7}).
\end{pf*}

\begin{remark} The results in Example~\ref{ex-1.1}, 3 and 4 come
from (\ref{7.32})
and an estimate of $\varphi$ as given in (\ref{1.23}).
\end{remark}

\begin{example}\label{exam-6.2}
The proper norm given in (\ref{bp.1q}) is useful in the study of
permanental fields of L\'evy processes because it requires that the
potential of the process, $u(x,y)$ is a function of $x-y$. The
following norms are proper norms that do not require this condition.
They are functions of the transition probability density, $p_{s}(x,y)$,
of a transient Markov process $X$ with reference measure $m$.
%
\begin{equation}
\|\nu\|_{w}:= \biggl( \int\int\biggl( \int w(x,y)w(y,z) \,d \nu(y)
\biggr)^{2} \,dm(x) \,dm(z) \biggr)^{1/2},\label{n.lj1q}
\end{equation}
where
%
\begin{equation}
w(x,y)= \int_{0}^{\infty} {p_{s}(x,y) \over\sqrt{\pi s}}
\,ds\label{n.lj2}
\end{equation}
and
%
\begin{equation}
\|\nu\|_{\Phi}:= \biggl( \int\int\Phi(x,y) \,d \nu(x) \,d \nu(y)
\biggr)^{1/2},\label{n1.9q}
\end{equation}
where $ \Phi(x,y)= \Theta_{l} (x,y) \Theta_{r} (x,y )$ and
\begin{eqnarray*}
\Theta_{l}(x,y)&=&\int_{0}^{\infty} \int
p_{s/2}(x,u)p_{s/2 }(y,u) \,dm(u) \,ds,
\\
\Theta_{r} (x,y)&=&\int_{0}^{\infty}
\int p_{s/2}(u,x )p_{s/2}(u,y) \,dm(u) \,ds.
\end{eqnarray*}

Proofs are given in an earlier version of this paper, with the same
title, \cite{LMR}, Section~6.
\end{example}
%
%

%



\printaddresses

\end{document}